\documentclass[reqno]{amsart}
\usepackage{amsmath, amssymb, eucal, amscd, amstext, enumerate}
\usepackage{amsfonts,amsthm}
\usepackage{mathrsfs}

\theoremstyle{plain}
\newtheorem{thm}{Theorem}[section]
\newtheorem{lem}[thm]{Lemma}
\newtheorem{eg}[thm]{Example}
\newtheorem{prop}[thm]{Proposition}
\newtheorem{cor}[thm]{Corollary}
\newtheorem{defn}[thm]{Definition}
\newtheorem{rem}[thm]{Remark}
\newtheorem{rem-ntn}[thm]{Remark and Notation}

\newtheorem{quest}[thm]{Question}

\newenvironment{prf}{{\noindent \textbf{Proof:}\ }}{\hfill $\Box$\\ \smallskip}

\numberwithin{equation}{section}

\newcommand{\smnoind}{{\smallskip\noindent}}

\newcommand{\ti}{\tilde}
\newcommand{\la}{\langle}
\newcommand{\ra}{\rangle}

\newcommand{\CL}{\mathcal{L}}
\newcommand{\CP}{\mathcal{P}}

\newcommand{\CU}{\mathcal{U}}

\newcommand{\CI}{\mathcal{I}}
\newcommand{\D}{\mathcal{D}}

\newcommand{\KH}{\mathfrak{H}}
\newcommand{\KK}{\mathfrak{K}}

\newcommand{\KI}{\mathfrak{I}}
\newcommand{\KJ}{\mathfrak{J}}

\newcommand{\Kl}{\mathfrak{l}}

\newcommand{\BC}{\mathbb{C}}
\newcommand{\BR}{\mathbb{R}}
\newcommand{\BN}{\mathbb{N}}

\newcommand{\BG}{{\mathbb{G}}}

\newcommand{\supp}{\mathrm{supp}\ \!}

\newcommand{\IT}{\mathfrak{T}}

\begin{document}

\title[]{A property for locally convex $^*$-algebras related to Property $(T)$ and character amenability}

\author{Xiao Chen \and Anthony To-Ming Lau \and Chi-Keung Ng}

\address[Xiao Chen]{Chern Institute of Mathematics, Nankai University, Tianjin 300071, China.}
\email{cxwhsdu@126.com}
\address[Anthony To-Ming Lau]{Department of Mathematicial and Statistical Sciences, University of Alberta, Edmonton, Alberta, Canada T6G-2G1}
\email{tlau@math.ualberta.ca}
\address[Chi-Keung Ng]{Chern Institute of Mathematics and LPMC, Nankai University, Tianjin 300071, China.}
\email{ckng@nankai.edu.cn; ckngmath@hotmail.com}


\begin{abstract}
For a locally convex $^*$-algebra $A$ equipped with a fixed continuous $^*$-character $\varepsilon$ (which is roughly speaking a generalized $F^*$-algebra), we define a cohomological property, called property $(FH)$, which is similar to character amenability.
Let $C_c(G)$ be the space of continuous functions on a second countable locally compact group $G$ with compact supports, equipped with the convolution $^*$-algebra structure and a certain inductive topology.
We show that $(C_c(G), \varepsilon_G)$ has property $(FH)$ if and only if $G$ has property $(T)$.
On the other hand, many Banach algebras equipped with canonical characters have property $(FH)$ (e.g., those defined by 
a nice locally compact quantum group). 
Furthermore, through our 
studies on both property $(FH)$ and character amenablility, we obtain 
characterizations of 
property $(T)$, amenability and compactness of $G$ in terms of the vanishing of one-sided cohomology of certain topological algebras, 
as well as in terms of fixed point properties. 
These three sets of characterizations can be regarded as analogues of one another. 
Moreover, we show that $G$ is compact if and only if the normed algebra $\big\{f\in C_c(G): \int_G f(t)dt =0\big\}$ (under $\|\cdot\|_{L^1(G)}$) admits a bounded approximate identity with the supports of all its elements being contained in a common compact set. 
\end{abstract}
\maketitle

\section{Introduction}

\medskip

The notion of property $(T)$ for locally compact groups was first introduced by Kazhdan in the 1960s (see \cite{Kah}) and was proved to be very useful.
A locally compact group $G$ is said to have \emph{property $(T)$} when every continuous unitary representation of $G$ having almost invariant unit vectors
actually has a non-zero invariant vector (see \cite[\S 1.1]{BHV}). 
Property $(T)$ has many equivalent formulations (see e.g.\ \cite{BHV}). 
In the case when the group is $\sigma$-compact, one equivalent form is given by the Delorme-Guichardet theorem (see \cite[Theorem 2.12.4]{BHV}): 
\begin{quotation}
a $\sigma$-compact locally compact group $G$ has property $(T)$ if and only if it has property $(FH)$, 
\end{quotation}
where property $(FH)$ can be viewed as the vanishing of the first cohomology $H^1(G, \pi)$ for any continuous unitary representation $\pi$ of $G$.

\medskip

In \cite{Kyed}, Kyed obtained a Delorme-Guichardet type theorem for a separable discrete quantum group $\BG$ as follows:  
if $({\rm Pol}(\widehat \BG), \Delta, S, \varepsilon)$ is the canonical Hopf $^*$-algebra associated with the dual compact quantum group $\widehat \BG$, then the vanishing of a certain cohomology of $({\rm Pol}(\widehat \BG), \varepsilon)$, with coefficients in $^*$-representations of ${\rm Pol}(\widehat \BG)$, is equivalent to the property $(T)$ of $\BG$ (as introduced by Fima in \cite{Fim}). 

\medskip

One may regard $({\rm Pol}(\widehat \BG), \varepsilon)$ as a ``non-commutative pointed space'', in the sense that 
it is a locally convex 
$^*$-algebra (in this case, the topology is the discrete one) equipped with a continuous $^*$-character. 
Motivated by 
this result of Kyed, we 
study certain cohomological properties for a general ``non-commutative pointed space'' $(A, \varepsilon)$
and partially generalize Kyed's result to the case of locally compact groups. 

\medskip

Note that there is a related notion of property $(T)$ for operator algebras (see e.g.\ \cite{CJ}  as well as \cite{Bek-T}, \cite{Bro}, \cite{LN}, \cite{LNW} and \cite{Ng-prop-T-gen}).
It was shown in \cite{Bek-T} that a discrete group has property $(T)$ if and only if its reduced group $C^*$-algebra has property $(T)$. 
This gives a characterization of property $(T)$ of a discrete group in terms of a certain $C^*$-algebra associated with it. 
In the case of a general locally compact group, if it has property $(T)$, then its reduced group $C^*$-algebra also has property $(T)$ (see \cite{Ng-prop-T-gen}).
However, the converse of this seems to be open. 

\medskip 

On the other hand, our quest for a partial generalization of Kyed's result
leads to a characterization of property $(T)$ of a second countable locally compact group $G$ in terms of what we called ``property $(FH)$'' of the locally convex $^*$-algebra $C_c(G)$ (under a certain inductive limit topology $\IT$), together with the canonical $^*$-character $\varepsilon_G$. 
Notice that property $(FH)$ can also be regarded as an analogue of character amenability for topological $^*$-algebras. 
We will also give some studies on both property $(FH)$ and character amenability. 
In particular, we investigate their relations with one-sided cohomology as well as with fixed point property.

\medskip
The following theorem collects some of the results of this paper. 
Note that the equivalence of Statements (T1) and (T4) as well as that of Statements (A1) and (A2) are well-known. 
We list them here for comparison with the other results in this paper. 
Note also that the equivalence of Statements (C1) and (C4) is implicitly established in the proof of Corollary \ref{cor:fix-pt} (and is not explicitly stated in the main contents of the paper).

\medskip

\begin{thm}\label{thm:summary}
Suppose that $G$ is a second countable locally compact group and $\varepsilon_G:C_c(G)\to \BC$ is the integration. 

\medskip\noindent
(a) The following statements are equivalent.
\begin{enumerate}[\ \ ({T}1).]
\item $G$ has property $(T)$. 

\item The non-commutative pointed space $(C_c(G), \varepsilon_G)$ has property $(FH)$.  

\item All the first left Hochschild cohomology of the locally convex $^*$-algebra 
\begin{equation}\label{eqt:defn-ti-C_c^0}
\ti C_c^0(G):=\left\{(f,\lambda)\in C_c(G)\oplus \BC: \lambda = -\int_G f(t)dt\right\}
\end{equation}
(under the topology induced by $\IT$) with coefficients in norm-continuous unitary representations of $G$ vanish. 

\item Every isometric affine-linear action of $G$ on a Hilbert space has a fixed point.
\end{enumerate}

\medskip\noindent
(b) The following statements are equivalent.
\begin{enumerate}[\ \ ({A}1).]
\item $G$ is amenable. 
	
\item The Banach algebra $(L^1(G), \|\cdot\|_{L^1(G)})$ is $\varepsilon_G$-amenable.
	
\item All the first left Hochschild cohomology of the Banach algebra 
\begin{equation}\label{eqt:defn-L^1_0}
\ti L^1_0(G):= \left\{(f,\lambda)\in L^1(G)\oplus \BC: \lambda = -\int_G f(t)dt\right\}
\end{equation}
with coefficients in ``dual Banach left $G$-modules'' vanish.

\item Any weak$^*$-continuous affine-linear action of $G$ on any dual Banach space with a norm-bounded orbit has a fixed point. 
\end{enumerate}

\medskip\noindent
(c) The following statements are equivalent.
\begin{enumerate}[\ \ ({C}1).]
\item $G$ is compact. 

\item The locally convex algebra $C_c(G)$ is $\varepsilon_G$-amenable. 

\item All the first left Hochschild cohomology of the locally convex algebra $\ti C_c^0(G)$
with coefficients in ``dual Banach left $G$-modules'' vanish.

\item Every continuous affine-linear action of $G$ on a quasi-complete locally convex space has a fixed point. 
\end{enumerate}
\end{thm}

\medskip

We also show that the amenability of $G$ is equivalent to the following statement. 
{\it
\begin{enumerate}[\ \ ({A}1').]
	\setcounter{enumi}{2}
	\item All the first left Hochschild cohomology of the Banach algebra 
	\begin{equation}\label{eqt:defn-ti-L^1_0}
	L^1_0(G):= \left\{f\in L^1(G): \int_G f(t)dt =0\right\}. 
	\end{equation} 
	with coefficients in ``dual Banach left $G$-modules'' vanish.
\end{enumerate}
}

\medskip

On the other hand, by \cite[Theorem 4.10]{Lau83}, $G$ is amenable if and only if there exists a bounded approximate identity in 
the Banach algebra $L^1_0(G)$. 
As an analogue of this fact, we show that the compactness of $G$ is equivalent to the following two statements concerning the 
algebra 
\begin{equation}\label{eqt:defn-C_c^0}
C_c^0(G):= C_c(G)\cap L^1_0(G). 
\end{equation}

{\it 
\begin{enumerate}[\ \ ({C}1).]
\setcounter{enumi}{4}
\item The locally convex subalgebra $C_c^0(G)$ of $C_c(G)$
has a bounded approximate identity.

\item The normed algebra $(C_c^0(G), \|\cdot\|_{L^1(G)})$
has a bounded approximate identity with the supports of all its elements being contained in a common compact set.
\end{enumerate}}

\bigskip

\section{Notations}

\medskip

Throughout this paper, all vector spaces, unless specified otherwise, are over the complex field (but most of the results have their counterparts in the field of real numbers). 
All topologies are Hausdorff, and all integrations, unless stated otherwise, are the Bochner integrations with respects to the norm topologies on the range spaces.

\medskip

If $S$ is a subset of a topological space $Y$, then we denote by $\overline{S}$ the closure of $S$ in $Y$.
For Banach spaces $E$ and $F$, we denote by $\CL(E;F)$ the set of all bounded linear operators from $E$ to $F$, and $\CL(E):=\CL(E;E)$. 

\medskip

By a \emph{locally convex algebra}, we mean an algebra $A$ equipped with a locally convex Hausdorff topology $\tau$ such that the multiplication is separately continuous.
We will denote this locally convex algebra by $A_\tau$. 
If, in addition, there is a continuous involution on $A$, then $A_\tau$ is called a \emph{locally convex $^*$-algebra}. 
Examples of locally convex $^*$-algebras include all Banach $^*$-algebras and the measure algebra for a locally compact group equipped with the weak$^*$-topology (when it is considered as the dual space of the space of continuous functions on the group that vanishes at infinity).
In the case when $A$ is a Banach $^*$-algebra with $\tau$ being the norm-topology, we may write $A$ instead of $A_\tau$. 

\medskip

Suppose that $\Phi$ is either a homomorphism or anti-homomorphism from 
$A$ to the algebra $L(Z)$ of linear maps on a vector space $Z$. 
For any subspace $X\subset A$ and $Y\subset Z$, we set $\Phi(X)Y$ to be the linear span of 
$\{\Phi(a)y:a\in X;y\in Y\}$.
Moreover, when $Z=A$, we denote $X\cdot Y :=\mathrm{m}(X)Y$ and $X^2 := \mathrm{m}(X)X$, where $\mathrm{m}:A\to L(A)$ is 
the homomorphism given by multiplication. 

\medskip

A subset $S\subseteq A$ is said to be \emph{$\tau$-bounded} if for any $\tau$-neighborhood $V$ of zero, there exists $\kappa >0$ such that 
$S\subseteq \kappa V$. 
A net $\{a_i\}_{i\in \KI}$ in $A$ is called a \emph{left} (respectively, \emph{right}) \emph{$\tau$-approximate identity} if 
$a_ib \stackrel{\tau}{\to} b$ (respectively, $ba_i \stackrel{\tau}{\to} b$), for any $b\in A$.  
Moreover, $\{a_i\}_{i\in \KI}$ is called an \emph{$\tau$-approximate identity} if it is both a left and a right $\tau$-approximate identity. 
A representation (respectively, an anti-representation) $\Psi:A\to \CL(E)$ of $A$ is said to be \emph{non-degenerate} if $\Psi(A)(E)$ is norm dense in $E$. 
A \emph{character} (respectively, \emph{$^*$-character}) on $A$ is a non-zero multiplicative linear (respectively, $^*$-linear) functional $\omega:A \to \BC$. 
A map from $A$ to a Banach space $F$ is said to be \emph{$\tau$-continuous} if it is continuous with respective to the norm-topology on $F$. 



\medskip

Furthermore, throughout this paper, $G$ is a locally compact group, $\Delta:G\to \mathbb{R}_+$ is the modular function and $\lambda^G$ is a fixed left Haar measure on $G$.
All integrations of maps on $G$ are with respect to $\lambda^G$ (although they are all written as $dt$ instead of $d\lambda^G(t)$). 
As usual, we use the convention that $\lambda^G(G) =1$ when $G$ is compact.

\medskip

For a topological space $(X, \CP)$, an action $\alpha$ of $G$ on $X$ is said to be \emph{$\CP$-continuous} if $\alpha_t:X\to X$ is $\CP$-$\CP$-continuous for each $t\in G$ and the map $t\mapsto \alpha_t(x)$ is $\CP$-continuous for any $x\in X$. 
A representation (or anti-representation) $\mu:G\to \CL(E)$ is said to be norm-continuous if 
the corresponding action on $E$ is norm-continuous. 

\medskip

We denote by $\KK(G)$ the collection of non-empty compact subsets of $G$.
For any function $f:G\to \BC$, we define $\supp f:= \overline{\{t\in G: f(t)\neq 0\}}$, and set
$$C_K(G)\ :=\ \{f:G\to \BC: f \text{ is continuous and } \supp f\subseteq K\} \qquad (K\in \KK(G)),$$
as well as $C_c(G):=\bigcup_{K\in \KK(G)} C_K(G)$. 
We consider $\varepsilon_G: C_c(G)\to \BC$ to be the $^*$-character induced by integration on $G$.

\medskip

Also, $L^1(G)$ is the Banach $^*$-algebra of (equivalence classes) of all $\lambda^G$-integrable complex functions, equipped with the $L^1$-norm $\|\cdot\|_{L^1(G)}$ as well as the convolution and the canonical involution:
$$(f\ast g)(t) := \int_G f(s)g(s^{-1}t)\ \!ds \quad \text{and}\quad f^*(t) := \Delta(t^{-1})\overline{f(t^{-1})}.$$ 
We regard $C_c(G)$ as a $^*$-subalgebra of $L^1(G)$. 
For each $f\in L^1(G)$, we define $\rho_s(f)(t) := \Delta(s)f(ts)$ and $\lambda_s(f)(t) := f(s^{-1}t)$ ($t\in G$).

\medskip

The convolution also gives a $^*$-representation of $L^1(G)$ on the Hilbert space $L^2(G)$ of square integrable functions, known as the ``left regular representation''.
The predual of the von Neumann subalgebra, $L(G)$, of $\CL(L^2(G))$ generated by the image of the left regular representation will be denoted by $A(G)$. 
Recall that $A(G)$ is the dense $^*$-subalgebra of the commutative $C^*$-algebra $C_0(G)$ consisting of coefficient functions of the left regular representation.
This $^*$-algebra structure on $A(G)$ turns it into a Banach $^*$-algebra under the predual norm, known as the \emph{Fourier algebra} of $G$ (see \cite{Eym} for more details).

\medskip

Let $\CU$ be an open neighborhood base of the identity $e\in G$ with each element being symmetric and having compact closure.
For each $U\in \CU$, we set
$$\CU(U):= \{V\in \CU: V\subseteq U\},$$
and fix a positive symmetric function $h_U$ with $\supp h_U\subseteq U$ and $\int_G h_U(t)\ dt = 1$.
It is well-known that $\{h_V\}_{V\in \CU}$ is a norm-bounded approximate identity for $L^1(G)$.
In fact, one has a slightly stronger property as stated in Lemma \ref{lem:equi-cont}.

\bigskip

\section{Topological $1$-cocycles of $C_c(G)$}

\medskip

In this section, we want to see when certain derivations of $C_c(G)$ come from continuous $1$-cocycles of $G$. 
This is related to the following analytic properties of the derivation. 

\medskip

\begin{defn}\label{defn:cont-loc-bdd}
For a Banach space $E$, a linear map $\theta: C_c(G)\to E$ is said to be:

\medskip
\noindent
(a) \emph{locally bounded}, if $\theta|_{C_K(G)}$ is $L^1$-bounded for any $K\in \KK(G)$.

\medskip
\noindent
(b) \emph{continuously locally bounded}, if there is a family $\{\kappa_K\}_{K\in \KK(G)}$ in $\BR_+$ satisfying
\begin{enumerate}[\ \ 1)]
\item $\|\theta(f)\|_E \leq \kappa_K\|f\|_{L^1(G)}$ ($K\in \KK(G); f\in C_K(G)$),
\item for any $\epsilon >0$, there is $U_\epsilon\in \CU$ such that $\kappa_K < \epsilon$ whenever $K\subseteq U_\epsilon$.
\end{enumerate}
\end{defn}

\medskip

Let $\mu:G\to \CL(E)$ be a norm-continuous representation (respectively, anti-representation) on a non-zero Banach space $E$. 
We denote by $\hat \mu:C_c(G)\to \CL(E)$ the representation (respectively, anti-representation) define by 
$$\hat \mu(f)(x) := \int_G f(t)\mu_t(x)\ \!dt \qquad (f\in C_c(G);x\in E).$$ 
Clearly, $\|\hat \mu(h_V)(x) - x\|\to 0$ ($x\in E$) and $\hat \mu$ is non-degenerate. 
Moreover, $\hat \mu$ is locally bounded. 
In fact, given $K\in \KK(G)$, since the set $\{\|\mu_t(x)\|:t\in K\}$  is bounded for any $x\in E$, the uniform boundedness principle produces $\kappa_K >0$ with $\|\mu_t\| \leq \kappa_K$ for any $t\in K$, and $\|\hat \mu(f)\| \leq \kappa_K \|f\|_{L^1(G)}$ whenever $f\in C_K(G)$. 
The same is true for anti-representation. 

\medskip

Conversely, it can be shown that any locally bounded non-degenerate representation (respectively, anti-representation) of $C_c(G)$ comes from a unique norm-continuous representation (respectively, anti-representation) of $G$ as above (see Lemma \ref{lem:cont-rep}).

\medskip

In the following, $\mu:G\to \CL(E)$ is a norm-continuous anti-representation.
We denote by $\pi:G\to \CL(E^*)$ the representation induced by $\mu$.
We also set $\Phi:= \hat\mu$ and $\check \Phi(f) := \Phi(f)^*\in \CL(E^*)$ ($f\in C_c(G)$). 
It is easy to see that $\check \Phi$ is a locally bounded representation.
By considering the approximate identity $\{h_V\}_{V\in \CU}$, the representation of $C_c(G)$ on the norm-closure of $\check \Phi(C_c(G))E^*$ induced by $\check \Phi$ is non-degenerate. 
On the other hand, when $E$ is reflexive, $\check \Phi$ is always non-degenerate, and 
hence $\pi$ is norm-continuous (by Lemma \ref{lem:cont-rep}).



\medskip

A map $c:G\to E^*$ is called a \emph{$1$-cocycle} for $\pi$ if 
\begin{equation}\label{eqt:gp-1-cocyc}
c(st)\ =\ \pi(s)(c(t)) + c(s)\qquad (s,t\in G), 
\end{equation}
and it is called a \emph{$1$-coboundary for $\pi$} if there is $\varphi\in E^*$ such that $c(t) = \pi(t)(\varphi) - \varphi$ ($t\in G$). 
A $1$-cocycle is said to be \emph{norm-continuous} (respectively, \emph{weak$^*$-continuous}) if it is continuous when $E^*$ is equipped with the norm-topology (respectively, the 
weak$^*$-topology).

\medskip

The following is the key lemma for the many results in this article, which gives a bijective correspondence between certain $1$-cocycles and certain ``derivations''. 

\medskip

\begin{lem}\label{lem:cf-1-cocyc}
Let $\mu$, $\pi$, $\Phi$ and $\check\Phi$ be as in the above.
Consider a weak$^*$-continuous $1$-cocycle $c:G\to E^*$ for $\pi$ and a locally bounded linear map $\theta:C_c(G)\to E^*$ satisfying
\begin{equation}\label{eqt:alg-cocyc}
\theta(f\ast g)\ =\ \check \Phi(f)\theta(g) + \theta(f)\varepsilon_G(g)
\qquad (f,g\in C_c(G)). 
\end{equation}

\smnoind
(a) The map $\CI(c): C_c(G) \to E^*$ given by 
\begin{eqnarray}\label{eqt:def-I-c}
\CI(c)(f)(x)\ :=\ \int_G f(t)c(t)(x)\ \!dt \qquad (f\in C_c(G); x\in E)
\end{eqnarray}
is locally bounded and satisfies \eqref{eqt:alg-cocyc}. 
If, in addition, $c$ is norm-continuous, then $\CI(c)$ is continuously locally bounded. 

\smnoind
(b) There is a weak$^*$-continuous $1$-cocycle $\D(\theta):G\to E^*$ for $\pi$ with $\theta = \CI(\D(\theta))$. 
If, in addition, $\theta$ is continuously locally bounded, then $\D(\theta)$ is norm-continuous.

\smnoind
(c) $\D(\CI(c)) = c$.

\smnoind
(d) $\theta$ is $L^1$-bounded if and only if $\{\D(\theta)(t): t\in G\}$ is a bounded subset of $E^*$.
\end{lem}
\begin{prf}
(a) Notice that if $K$ is a non-empty compact subset of $G$ and $x\in E$, then the set $\{c(t)(x):t\in K\}$ is bounded, and the uniform boundedness principle tells us that 
\begin{equation}\label{eqt:weak-star-cont>loc-bdd}
\kappa^c_K\ := \ {\sup}_{t\in K} \|c(t)\|\ <\ \infty. 
\end{equation}
Thus, $\|\CI(c)(f)\| \leq \kappa^c_K\|f\|_{L^1(G)}$ for any $f\in C_K(G)$, and $\CI(c)$ is locally bounded.  

Moreover, Equality \eqref{eqt:gp-1-cocyc} gives
\begin{eqnarray*}
\CI(c)(f\ast g)
& = & \int_G \int_G f(s)g(r)c(sr)\ \!dr\ \!ds\\
& = & \int_G \int_G f(s)\pi(s)g(r)c(r)\ \!dr\ \!ds + \int_G \int_G f(s)g(r)c(s)\ \!dr\ \!ds\\
& = & \int_G f(s)\pi(s)\CI(c)(g)\ \!ds + \int_G f(s)c(s)\varepsilon_G(g)\ \!ds,
\end{eqnarray*} 
which means that $\CI(c)$ satisfies Relation \eqref{eqt:alg-cocyc}.

Now, suppose that $c$ is norm-continuous. 
Since $c(e) = 0$, for any $\epsilon >0$, one can find $U\in \CU$ such that $\kappa_K^c\leq \epsilon$ whenever $K\subseteq U$, and $\CI(c)$ is continuously locally bounded.

\smnoind
(b) Fix $f\in C_c(G)$, $s\in G$  and $U\in \CU$.
Let $W$ be an open neighborhood of $s$ with compact closure and set $L:= (\supp f) \cdot \overline{W} \cdot \overline{U}$.
Consider the net $\{h_V\}_{V\in \CU(U)}$ as in Section 2.
Then both $\rho_{r^{-1}}(f)$ and $\rho_{r^{-1}}(f)\ast h_V = f\ast\lambda_r(h_V)$ belongs to $C_L(G)$ when $r\in W$ and $V\in \CU(U)$.
The local boundedness of $\theta$ means that one can find $\kappa > 0$ with $\|\theta(g)\|_{E^*}\leq \kappa \|g\|_{L^1(G)}$ for every $g\in C_L(G)$.
Thus,
\begin{equation}\label{eqt:unif-conv-2}
{\sup}_{r\in W} \|\theta(f\ast \lambda_r(h_V)) - \theta(\rho_{r^{-1}}(f))\| \leq  \kappa \cdot {\sup}_{r\in W}\|\rho_{r^{-1}}(f)\ast h_V - \rho_{r^{-1}}(f)\|.
\end{equation}

For any $n\in \BN$, $f_1,...,f_n\in C_c(G)$, $x_1,...,x_n\in E$ and $r\in W$, Relation \eqref{eqt:alg-cocyc} gives
\begin{eqnarray*}
\theta\big(\lambda_r(h_V)\big) \Big({\sum}_{k=1}^n \Phi(f_k)(x_k)\Big)
& = & {\sum}_{k=1}^n  \check \Phi(f_k) \theta\big(\lambda_r(h_V)\big)(x_k) \\
& = & {\sum}_{k=1}^n  \big(\theta\big(f_k\ast \lambda_r(h_V)\big) - \theta(f_k)\varepsilon_G\big(\lambda_r(h_V)\big)\big)  (x_k),
\end{eqnarray*}
which converges uniformly in the varaible $r\in W$ to ${\sum}_{k=1}^n  \big( \theta(\rho_{r^{-1}}(f_k)) - \theta(f_k)\big)(x_k) $, because of \eqref{eqt:unif-conv-2} and Lemma \ref{lem:equi-cont}.
Moreover, one has $\big\|\theta\big(\lambda_r(h_V)\big)\big\|_{E^*}\leq \kappa$ whenever $r\in W$ and $V\in \CU(U)$.
Hence, the norm density of $\Phi(C_c(G))E$ in $E$ will imply that for each $x\in E$, the net
$\big\{\theta(\lambda_r(h_V))(x)\big\}_{V\in \CU(U)}$ is uniformly Cauchy for all $r\in W$.

Consequently, for $r\in W$, we may define a linear map $\D(\theta)(r):E\to \BC$ by 
\begin{equation}\label{eqt:defn-D-theta}
\D(\theta)(r)(x)\ :=\ {\lim}_V \theta(\lambda_r(h_V))(x)\qquad (x\in E).
\end{equation}
In particular,
\begin{equation}\label{eqt:dense-subsp-2}
\D(\theta)(r)\Big({\sum}_{k=1}^n \Phi(f_k)(x_k)\Big)
\ =\ {\sum}_{k=1}^n  \big(\theta\big(\rho_{r^{-1}}(f_k)\big) - \theta(f_k)\big)(x_k).
\end{equation}
Since $\|\D(\theta)(r)\| \leq \kappa$ (because of \eqref{eqt:defn-D-theta}), one knows that $\D(\theta)(r)\in E^*$, and that $\D(\theta)(r)$ does not depend on the choices of $U$, $\{h_V\}_{V\in\CU(U)}$ nor $W$ (so long as $W$ contains $r$).

Furthermore, as the function $r\mapsto \theta\big(\lambda_r(h_V)\big)(x)$ is continuous on $W$ for all $V\in \CU(U)$, and $\big\{\theta(\lambda_r(h_V))(x)\big\}_{V\in \CU(U)}$ converges uniformly to $\D(\theta)(r)(x)$ for all $r\in W$, we see that $r\mapsto \D(\theta)(r)(x)$ is a continuous complex function on $W$, for any $x\in E$.
Consequently, $\D(\theta): G\to E^*$ is weak$^*$-continuous.

It is easy to check, by using \eqref{eqt:dense-subsp-2}, that
$\D(\theta)(st)(x) = \big(\pi(s)\D(\theta)(t)+\D(\theta)(s)\big) (x)$ whenever $s,t\in G$ and $x\in \Phi(C_c(G))E$.
Hence, $\D(\theta)$ satisfies \eqref{eqt:gp-1-cocyc} (as $\Phi$ is non-degenerate).

In order to verify $\theta = \CI(\D(\theta))$, we observe from \eqref{eqt:alg-cocyc} and \eqref{eqt:dense-subsp-2} that if $g\in C_c(G)$, 
\begin{eqnarray*}
\theta(g)\big(\Phi(f)(x)\big)
& = & \theta\Big(\int_G  f(r)\lambda_r(g) \ \!dr\Big)(x) - \theta(f)\varepsilon_G(g)(x)\\
& = & \int_G \theta\big(f(r)\lambda_r(g)\big)(x)\ \!dr - \int_G  g(s)\theta(f)(x) \ \!ds\\
& = & \int_G g(s)  \big(\theta(\rho_{s^{-1}}(f)) - \theta(f)\big)(x) \ \!ds\\
& = & \int_G g(s)  \D(\theta)(s) \big(\Phi(f)(x) \big) \ \!ds
\end{eqnarray*}
for any $f\in C_c(G)$ and $x\in E$ 
(the second equality above follows from the $L^1$-boundedness of $\theta|_{C_{K_f\cdot K_g}(G)}$ where $K_f := \supp f$ and $K_g:= \supp g$; notice that $\supp f(r)\lambda_r(g) \subseteq K_f\cdot K_g$, for any $r\in K_f$).
Again, 	the non-degeneracy of $\Phi$ gives the required equality. 

Finally, suppose that $\theta$ is continuously locally bounded and consider $\{\kappa_K\}_{K\in \KK(G)}$ to be the family  as in Definition \ref{defn:cont-loc-bdd}(b).
Let $\{s_j\}_{j\in \KJ}$ be a net in $G$ that converges to $e$.
For any $\epsilon >0$, let $U_\epsilon$ be as in Definition \ref{defn:cont-loc-bdd}(b)(2). 
Pick any $V_\epsilon\in \CU$ with $\overline{V_\epsilon} \cdot \overline{V_\epsilon}\subseteq U_\epsilon$.
There exists $j_0\in \KJ$ such that $s_j\in V_\epsilon$ whenever $j\geq j_0$.  
In the argument above, if we put $s=e$ and set the neighborhoods $U$ and $W$ to be $V_\epsilon$, then for any $V\in \CU(U)$ and $j\geq j_0$, one has $\supp \lambda_{s_j}(h_V) \subseteq U_\epsilon$, which implies $\|\theta(\lambda_{s_j}(h_V))\| < \epsilon$.  
Thus, Relation \eqref{eqt:defn-D-theta} tells us that $\|\D(\theta)(s_j)\|  \leq \epsilon$ and $\D(\theta)$ is continuous at $e$.
Now, if $\{t_j\}_{j\in \KJ}$ is a net in $G$ converging to $t\in G$, then 
$\D(\theta)(t_j) = \pi(t_jt^{-1})\D(\theta)(t) + \D(\theta)(t_jt^{-1}) \to \D(\theta)(t)$. 

\smnoind
(c) For any $r\in G$, $f\in C_c(G)$ and $x\in E$, one has, by Equalities \eqref{eqt:dense-subsp-2} and \eqref{eqt:gp-1-cocyc}, 
\begin{eqnarray*}
\D\big(\CI(c)\big)(r)\big(\Phi(f)x\big)
& = & \Big(\CI(c)\big(\rho_{r^{-1}}(f)\big) - \CI(c)(f)\Big) (x) \\
& = & \int_G \big(\rho_{r^{-1}}(f)(s)c(s)\big)(x)\ \!ds - \int_G f(t)c(t)(x)\ \!dt\\
& = & \int_G f(t)c(tr)(x)\ \!dt - \int_G f(t)c(t)(x)\ \!dt\\
& = & \int_G f(t)\pi(t)(c(r))(x)\ \!dt
\quad = \quad c(r)\big(\Phi(f)(x)\big).
\end{eqnarray*}
As $\Phi$ is non-degenerate, we know that $\D(\CI(c))(r) = c(r)$. 

\smnoind
(d) If $\theta$ is $L^1$-bounded, then Relation \eqref{eqt:defn-D-theta} implies that $\|\D(\theta)(t)(x)\| \leq \|\theta\| \|x\|$ ($x\in E$), i.e.\ $\|\D(\theta)(t)\|\leq \|\theta\|$ ($t\in G$).
The converse follows directly from \eqref{eqt:def-I-c} and the equality $\theta = \CI(\D(\theta))$.
\end{prf}

\medskip

Suppose that $\KH$ is a Hilbert space and $\pi:G\to \CL(\KH)$ is a norm-continuous unitary representation.
Then $\pi$ is induced by a norm-continuous anti-representation $\mu$ of $G$ on $E:=\KH^*$. 
It is well-known that  a weakly measurable $1$-cocycle for $\pi$ is automatically continuous when $\KH$ is separable (see e.g.\ \cite[2.14.3]{BHV}).
Moreover, every bounded $1$-cocycle for $\pi$ is a $1$-coboundary as $E^*$ is a Hilbert space (see \cite[Proposition 2.2.9]{BHV}).
These, together with Lemma \ref{lem:cf-1-cocyc} and the Delorme-Guichardet theorem, produce the following result.

\medskip

\begin{thm}\label{thm:alter-prop-T-1}
(a) If $G$ has property $(T)$ (respectively, and is second countable), then for any norm-continuous unitary representation $\pi$ of $G$ on a Hilbert space $\KH$, any continuously locally bounded map (respectively, any locally bounded map) $\theta:C_c(G)\to \KH$ satisfying \eqref{eqt:alg-cocyc} is $L^1$-bounded.

\smnoind
(b) If $G$ is $\sigma$-compact such that for any norm-continuous unitary representation $\pi$ of $G$ on a Hilbert space $\KH$, every continuously locally bounded map $\theta:C_c(G)\to \KH$  satisfying \eqref{eqt:alg-cocyc} is $L^1$-bounded, then $G$ has property $(T)$. 
\end{thm}

\medskip

In particular, \emph{when $G$ is second countable, $G$ has property $(T)$ if and only if for any norm-continuous unitary representation $\pi$ of $G$ on a Hilbert space $\KH$, every locally bounded map $\theta:C_c(G)\to \KH$ satisfying \eqref{eqt:alg-cocyc} is automatically $L^1$-bounded.}
One may restate this statement in terms of a certain cohomology theory for  locally convex $^*$-algebras with fixed $^*$-characters, as introduced in the following section.

\bigskip

\section{Property $(FH)$ and character amenability}

\medskip

In this section, we will introduce and study some cohomological properties for locally convex $^*$-algebras. 
In particular, we will establish their relations with character amenability (for topological algebras) as well as with property $(T)$ (for locally compact groups). 

\medskip

\emph{Throughout this section, $B_\sigma$ is a locally convex algebra, $A_\tau$ is a locally convex $^*$-algebra, $\omega:B\to \BC$ is a $\sigma$-continuous character and $\varepsilon:A\to \BC$ is a $\tau$-continuous $^*$-character.}

\medskip

For a non-zero Banach space $E$ and an anti-representation $\Phi:B\to \CL(E)$, we denote by $\check\Phi: B\to \CL(E^*)$ the homomorphism induced by $\Phi$. 
Observe that if $\Phi$ is $\sigma$-continuous, then $\check \Phi$ is also $\sigma$-continuous (since $\|\check \Phi(b)\| \leq \|\Phi(b)\|$). 

\medskip

Let us first extend the definition of character amenability, as studied in \cite{KLP08} and \cite{Mon08}, to the setting of locally convex algebras as follows.

\medskip

\begin{defn}\label{defn:char-amen}
(a) When $\Phi:B \to \CL(E)$ and $\Psi:B \to \CL(F)$ are representations of $B$ on Banach spaces, a linear map $\theta: B\to \CL(F;E)$ is called a \emph{$1$-cocycle for $(\Phi, \Psi)$} (or a \emph{$(\Phi, \Psi)$-derivation}) if 
\begin{equation}\label{eqt:deriv}
\theta(ab)\ =\ \Phi(a)\theta(b) + \theta(a)\Psi(b) \qquad (a,b\in B), 
\end{equation}
and it is called a \emph{$1$-coboundary}  (or is said to be \emph{inner}) if 
there is $T\in \CL(F;E)$ with 
\begin{equation}\label{eqt:inner}
\theta(a)\ =\ T\Psi(a) - \Phi(a)T \qquad (a\in B).
\end{equation}

\smnoind
(b) $B_\sigma$ is said to be \emph{$\omega$-amenable} if for each $\sigma$-continuous anti-representation $\Phi:B\to \CL(E)$, every $\sigma$-continuous $(\check\Phi,\omega)$-derivation $\theta:B\to E^*$ is inner.

\smnoind
(c) $B_\sigma$ is said to have \emph{property $(CB)$} if for any $\sigma$-continuous anti-representation $\Phi:B\to \CL(E)$ and any $\sigma$-continuous 
left $B$-module map $\Xi:B\to E^*$, one can find $\varphi_0\in E^*$ with $\Xi(a) = \check \Phi(a)(\varphi_0)$ ($a\in B$). 
\end{defn}

\medskip

In part (c) above, we consider $E^*$ as a left $B$-module through the representation of $B$ on $E^*$ induced by $\Phi$. 
Similarly, we regard $\KH$ as a left $A$-module in Definition \ref{defn:FH}(b) below, concerning a weaker cohomological property for locally convex 
$^*$-algebras. 

\medskip

\begin{defn}\label{defn:FH}
(a) We say that $(A_\tau, \varepsilon)$ has \emph{property $(FH)$} if for any $\tau$-continuous $^*$-representation 
$\Phi$ of $A$ on a Hilbert space, any $\tau$-continuous  $(\Phi, \varepsilon)$-derivation is inner.
	
\smnoind
(b) We say that $A_\tau$ satisfies \emph{property $(CH)$} if for any $\tau$-continuous $^*$-representation $\Phi:A\to \CL(\KH)$ and any $\tau$-continuous 
left $A$-module map  $\Lambda:A\to \KH$, one can find $\xi_0\in \KH$ such that $\Lambda(a) = \Phi(a)\xi_0$ ($a\in A$).
\end{defn}

\medskip

\begin{rem}\label{rem:sep}
It is well-known that character amenability and property $(FH)$ can be expressed in terms of the vanishing of a certain first Hochschild cohomology. 
	In the same way, properties $(CH)$ and $(CB)$ can be expressed in terms of the vanishing of the first ``left Hochschild cohomology''. 
	
	More precisely, suppose that $F$ is a Banach space and $\Psi: B\to \CL(F)$ is a $\sigma$-continuous representation. 
	One may consider the space $C^n(B;F)$ of separately continuous $n$-linear maps from the $n$-times product of $B$ to $F$ and define $\partial_n:C^n(B;F) \to C^{n+1}(B;F)$ by 
	$$\partial_n \Xi (b_0,...,b_n)\ :=\ \Psi(b_0)\Xi(b_1,...,b_n) - \Xi(b_0b_1,...,b_n) + \cdots + (-1)^n\Xi(b_0,b_1,...,b_{n-1}b_n).$$
	The vector space $H^n(B;\Psi):= \ker \partial_n / {\rm Im}\ \!\partial_{n-1}$ is called the \emph{$n$-th left Hochschild cohomology} with coefficients in $\Psi$. 
	Obviously, $B_\sigma$ has property $(CB)$ if and only if $H^1(B;\check\Phi) = (0)$ for any $\sigma$-continuous anti-representation $\Phi$ of $B$.
	
	Similarly, $A_\tau$ has property $(CH)$ if and only if $H^1(A;\Phi) = (0)$ for any $\tau$-continuous $^*$-representation $\Phi:A\to \CL(\KH)$. 
\end{rem}

\medskip

One may wonder whether it is possible to consider only non-degenerate anti-representations in the study of $\omega$-amenability. 
The following lemma tells us that this can be done when $B$ has a bounded approximate identity. 

\medskip

\begin{lem}\label{lem:char-amen-bai}
Suppose that $B$ has a $\sigma$-bounded $\sigma$-approximate identity $\{e_i\}_{i\in \KI}$.
Then $B_\sigma$ is $\omega$-amenable if and only if for every $\sigma$-continuous non-degenerate anti-representation $\Psi:B\to \CL(F)$, each $\sigma$-continuous $(\check\Psi,\omega)$-derivation is inner. 
\end{lem}
\begin{prf}
It suffices to show that the condition concerning non-degenerate anti-representations implies the $\omega$-amenability of $B_\sigma$. 
Suppose that $\Phi:B \to \CL(E)$ is a $\sigma$-continuous anti-representation and $\theta:B\to E^*$ is a $\sigma$-continuous $(\check\Phi, \omega)$-derivation.
\emph{In the following, we consider $a$, $x$ and $\varphi$ to be arbitrary elements in $B$, $E$ and $E^*$, respectively.}

First of all, we set $F$ to be the closure of $\Phi(B)E$. 
Clearly, $F$ is $\Phi$-invariant. 
As $\{e_i\}_{i\in \KI}$ is a $\sigma$-approximate identity, we know that $\Phi(B)E$ is a subset of the closure of $\Phi(B)F$, which implies that the induced anti-representation $\Psi:B\to \CL(F)$ is non-degenerate. 
Now, consider $P: E^* \to F^*$ to be the canonical map given by restrictions. 

Observe that 
$$\check \Phi(a)(\varphi)(u) = \varphi\big(\Phi(a)u\big) = P(\varphi)\big(\Psi(a)u\big) = \check\Psi(a)(P(\varphi))(u) \qquad (u\in F).$$ 
From this, one can check easily that $P\circ \theta$ is a $(\check\Psi, \omega)$-derivation, and the hypothesis produces $\psi_0\in F^*$ satisfying
\begin{equation}\label{eqt:P-Theta-inner}
\theta(a)(u) = P\circ \theta(a)(u) = \omega(a) \psi_0(u) - \check{\Psi}(a)(\psi_0)(u) \qquad (u\in F).
\end{equation}

Obviously, $\omega(e_i)\to 1$ (as there exists $b\in B$ with $\omega(b) \neq 0$). 
Moreover, since $\{e_i\}_{i\in \KI}$ is $\tau$-bounded and both $\theta$ and $\Phi$ are $\tau$-continuous, by considering a subnet if necessary, we may assume that $\{\theta(e_i)\}_{i\in \KI}$ weak$^*$-converges to an element $\chi_0\in E^*$ and $\{\check \Phi(e_i)\}_{i\in \KI}$ weak$^*$-converges to an operator $T_0\in \CL(E^*) = (E^*\otimes_\pi E)^*$. 
Then 
\begin{equation}\label{eqt:h>F-bot}
\chi_0(\Phi(a)x) = \lim_i \check \Phi(a)(\theta(e_i)) (x) = 
\lim_i \theta(ae_i)(x) - \theta(a)(x)\omega(e_i) = 0,
\end{equation}
and 
\begin{equation}\label{eqt:T-check-Phi-a-x}
T_0\big(\check \Phi(a)(\varphi)\big)(x) = \lim_i \check \Phi(e_i)\big(\check \Phi(a)(\varphi)\big)(x) = \check \Phi(a)(\varphi)(x).
\end{equation}
Furthermore, 
if we consider $E\subseteq E^{**}$ in the canonical way, then 
$$T_0^*(x)(\psi) = \lim_i \psi\big(\Phi(e_i)x\big) = 0 \quad \text{whenever} \quad \psi\in F^\bot,$$
where $F^\bot:=\{\psi\in E^*: \psi(F) = \{0\}\}$. 
This shows that $T_0^*(x) \in (F^\bot)^\bot\subseteq E^{**}$ and there is a net $\{u_j\}_{j\in \KJ}$ in $F$ such that 
\begin{equation}\label{eqt:T*-x-lim}
T_0(\varphi)(x) = T_0^*(x)(\varphi) = \lim_j \varphi(u_j).
\end{equation}
 
Now, consider any two elements $\psi_0', \psi_0''\in E^*$ satisfying $\psi_0'|_{F} = \psi_0''|_{F} = \psi_0$. 
By \eqref{eqt:T*-x-lim}, \eqref{eqt:P-Theta-inner} and \eqref{eqt:T-check-Phi-a-x}, 
\begin{eqnarray*}
T_0(\theta(a))(x) 
& = & \lim_j \theta(a)(u_j)
\ =\ \lim_j \omega(a)\psi_0(u_j) - \check \Psi(a)(\psi_0)(u_j) \\
& = & \omega(a) T_0\big(\psi_0'\big)(x) - \lim_j \psi_0''\big(\Phi(a) u_j\big) \\
& = & \omega(a) T_0\big(\psi_0'\big)(x) - T_0\big(\check \Phi(a)(\psi_0'')\big)(x)\\
& = & \omega(a) T_0\big(\psi_0'\big)(x) - \check \Phi(a)(\psi_0'')(x)
\end{eqnarray*}
This implies that
\begin{eqnarray*}
\theta(a)(x) 
& = & \lim_i \theta(e_ia)(x)
\ = \ \lim_i \check \Phi(e_i)\big(\theta(a))(x) + \omega(a)\theta(e_i)(x)\\
& = & T_0(\theta(a))(x) + \omega(a) \chi_0(x) \\
& = & \omega(a)\big(T_0(\psi_0')+\chi_0\big)(x) - \check \Phi(a)(\psi_0'')(x). 
\end{eqnarray*}

On the other hand, for each $u\in F$, one has
$$T_0(\psi_0')(u) = \lim_i \check \Phi(e_i)(\psi_0')(u) = \lim_i \psi_0'\big(\Phi(e_i) u\big) = \psi_0(u),$$
and 
$\chi_0(u) = 0$ because of \eqref{eqt:h>F-bot}. 
This shows that 
$T_0(\psi_0')+\chi_0\in E^*$ is an extension of $\psi_0$. 

Consequently, if we fix any extension 
$\psi_0'\in E^*$ of $\psi_0$ and set $\psi_0''= T_0(\psi_0')+\chi_0$, then we have 
$\theta(a) = \omega(a)\psi_0'' - \check \Phi(a)(\psi_0'')$ as required. 
\end{prf}

\medskip

We have an analogue of the above in the case of property $(FH)$ as stated in part (b) of the following lemma.

\medskip

\begin{lem}\label{lem:sep-n.d.}
	(a) If $A_\tau$ is separable, then $(A_\tau,\varepsilon)$ has property $(FH)$ if and only if for any $\tau$-continuous $^*$-representation $\Phi$ of $A$ on a separable Hilbert space, any $\tau$-continuous  $(\Phi, \varepsilon)$-derivation is inner.
	
	\smnoind
	(b) If $A$ has a $\tau$-bounded (left) $\tau$-approximate identity $\{e_i\}_{i\in \KI}$, then $(A_\tau, \varepsilon)$ has property $(FH)$ if and only if for any $\tau$-continuous non-degenerate $^*$-representation 
	$\Psi$ of $A$, any $\tau$-continuous $(\Psi, \varepsilon)$-derivation is inner.
	\end{lem}
\begin{prf}
(a) We only need to establish the sufficiency. 
Let $\Phi$ be a $\tau$-continuous $^*$-representation of $A$ on a Hilbert space $\KH$ (which may be inseparable) and $\theta$ be a $\tau$-continuous $(\Phi, \varepsilon)$-derivation.
The closure, $\KH_\theta$, of $\theta(A)$ is a separable $\Phi$-invariant subspace of $\KH$ (as $\theta$ satisfies \eqref{eqt:deriv}), and we let $\Phi_\theta : A\to \CL(\KH_\theta)$ be the induced $\tau$-continuous $^*$-representation. 
Obviously, if $\theta$ is inner as a $(\Phi_\theta, \varepsilon)$-derivation, then it is inner as a $(\Phi, \varepsilon)$-derivation, and part (a) is established. 
	
\smnoind
(b) It suffices to show that the condition concerning non-degenerate $^*$-representation implies property $(FH)$. 
In fact, suppose that $\Phi:A \to \CL(\KH)$ is a $\tau$-continuous $^*$-representation and $\theta$ is a $\tau$-continuous $(\Phi, \varepsilon)$-derivation.
We set $\KH_\Phi$ to be the closure of $\Phi(A)\KH$. 
Clearly, $\KH_\Phi$ is $\Phi$-invariant. 
By considering $\{e_i\}_{i\in \KI}$, we know that $\Phi(A)\KH$ is a subset of the closure of $\Phi(A)\KH_\Phi$, and hence the induced $^*$-representation $\Phi_0:A\to \CL(\KH_\Phi)$ is non-degenerate. 
If $P:\KH\to \KH_\Phi$ is the orthogonal projection, then $P\circ \theta$ is a $(\Phi_0,\varepsilon)$-derivation and hence there is $\xi_0\in \KH_\Phi$ 
satisfying 
$$P\circ \theta(a) = \varepsilon(a)\xi_0 - \Phi_0(a)\xi_0 \qquad (a\in A). $$
On the other hand, there is a subnet of $(1-P)(\theta(e_i))$ that weakly converges to some $\xi_1\in (1-P)(\KH)$.
Since $\Phi(A)(1-P)(\KH) = \{0\}$ 
and 
$$(1-P)\circ \theta (ab) = (1-P)\circ \theta(a)\varepsilon(b) \qquad (a,b\in A),$$ 
it is 
easy to check that $\theta(a) = \varepsilon(a)(\xi_0 + \xi_1) - \Phi(a)(\xi_0 + \xi_1)$ ($a\in A$).
\end{prf}

\medskip

It follows from Lemma \ref{lem:bai-B-psi}(b) that if $A_\tau$ has a $\tau$-bounded left (or right) $\tau$-approximate identity, then it has a $\tau$-bounded $\tau$-approximate identity (since $A$ has a $\tau$-continuous involution). 

\medskip

\begin{rem}
If $A_\tau$ is separable and has a $\tau$-bounded $\tau$-approximate identity, then it follows from Lemma \ref{lem:sep-n.d.}(a) and the argument of Lemma \ref{lem:sep-n.d.}(b) that $A_\tau$ has property $(FH)$ if and only if 
for any $\tau$-continuous non-degenerate $^*$-representation $\Phi$ of $A$ on a separable Hilbert space, any $\tau$-continuous  $(\Phi, \varepsilon)$-derivation is inner.
\end{rem}

\medskip

\begin{lem}\label{lem:two-elem-inner}
(a) Suppose that $\Phi:B\to \CL(E)$ is a non-degenerate anti-representation and $\theta: B\to E^*$ is a $(\check \Phi,\omega)$-derivation. 
If there exist $\varphi_1, \varphi_2\in E^*$ with $\theta(a) = \varphi_1\omega(a) - \check\Phi(a)(\varphi_2)$ ($a\in B$), then $\varphi_1 = \varphi_2$ (and hence $\theta$ is inner).
	
\smnoind
(b) If $B$ has a $\sigma$-bounded right $\sigma$-approximate identity 
$\{b_j\}_{j\in \KJ}$, then $B_\sigma$ has property $(CB)$. 
\end{lem}
\begin{prf}
(a) Equality \eqref{eqt:deriv} tells us that for any $a,b\in B$, one has 
	$$\varphi_1\omega(ab) - \check \Phi(ab)(\varphi_2) 
	\ = \ \check \Phi(a)\big(\varphi_1\omega(b) - \check\Phi(b)(\varphi_2)\big) + 
	\big(\varphi_1\omega(a) - \check \Phi(a)(\varphi_2)\big) \omega(b)$$
and hence  $\check \Phi(a)(\varphi_1) = \check \Phi(a)(\varphi_2)$ (as $\omega$ is non-zero).
Since $\Phi$ is non-degenerate, we conclude that $\varphi_1 = \varphi_2$.
	
\smnoind
(b) Suppose that $\Phi:B\to \CL(E)$ is a $\sigma$-continuous anti-representation of $B$ and $\Xi:B\to E^*$ is a $\sigma$-continuous left $B$-module map. 
As $\{\Xi(b_j)\}_{j\in \KJ}$ is norm-bounded, one can find a subnet $\{b_{j_i}\}_{i\in \KI}$ such that $\Xi(b_{j_i}) \to \varphi_0\in E^*$ under the weak$^*$-topology.  
For any $a\in B$, the net $\{\Xi (ab_{j_i})\}_{i\in \KI}$ will weak$^*$-converge to both $\Xi(a)$ and $\check \Phi(a)(\varphi_0)$. 
This gives the required conclusion. 
\end{prf}

\medskip

In the following, we set $\ti B := B\oplus \BC$ to be the unitalization of $B$ 
(whether or not $B$ is unital) and consider $\ti \sigma$ to be the direct sum topology on $\ti B$.
For any $\sigma$-continuous anti-representation $\Phi:B\to \CL(E)$, we denote by $\ti \Phi:\ti B\to \CL(E)$ the unital $\ti\sigma$-continuous anti-representation extending $\Phi$. 

\medskip

\begin{prop}\label{prop:CB<->char-amen}
Let $B^\omega := \ker \omega$.

\smnoind
(a) If $B$ has a $\sigma$-bounded $\sigma$-approximate identity and $B^\omega_{\sigma}$ has property (CB), then $B_\sigma$ is $\omega$-amenable. 

\smnoind
(b) If $B$ is unital and $B_\sigma$ is $\omega$-amenable,  then $B^\omega_{\sigma}$ has property $(CB)$ and $(B^\omega)^2$ is $\sigma$-dense in $B^\omega$. 

\smnoind
(c) The following statements are equivalent. 
\begin{enumerate}
\item $B_\sigma$ is $\omega$-amenable. 

\item $\ti B_{\ti \sigma}$ is $\ti \omega$-amenable. 

\item $\ti B^{\ti \omega}_{\ti \sigma}$ has property $(CB)$. 
\end{enumerate}
\end{prop}
\begin{prf}
(a) Let $\Phi:B\to \CL(E)$ be a $\sigma$-continuous non-degenerate anti-representation and $\theta:B\to E^*$ be a $\sigma$-continuous $(\check\Phi, \omega)$-derivation. 
Set $\Psi := \Phi|_{B^\omega}: B^\omega\to \CL(E)$ and $\Xi:= \theta|_{B^\omega}$. 
Then $\Xi$ is a left $B$-module map and the hypothesis produces $f_0\in E^*$ such that $\Xi(x) = \check\Psi(x)(f_0)$ ($x\in B^\omega$).
Pick any $u\in B$ with $\omega(u) = 1$. 
For any $x\in B^\omega$ and $\lambda\in \BC$, one has 
$$\theta(x+\lambda u)
\ = \ \Xi(x) + \lambda \theta(u) 
\ = \ \check \Phi(x+\lambda u)(f_0) - \big(\check \Phi(u)(f_0) - \theta(u)\big)\omega(x+\lambda u)$$
and Lemma \ref{lem:two-elem-inner}(a) implies that $\theta$ is inner. 
Now, Lemma \ref{lem:char-amen-bai} gives the conclusion. 

\smnoind
(b) As $B$ is unital, we have $B_\sigma = \widetilde{ B^\omega}_{\tilde{\sigma}}$. 
Let $\Psi:B^\omega\to \CL(E)$ be a $\sigma$-continuous anti-representation and $\Xi:B^\omega\to E^*$ be a $\sigma$-continuous left $B^\omega$-module map. 
Then $\Phi:=\ti \Psi$ is a unital $\sigma$-continuous anti-representation of $B$ on $E$. 
If $\theta:B\to E^*$ is defined by $\theta(x+\alpha 1) := \Xi(x)$ ($x\in B^\omega;\alpha\in \BC$), then $\theta$ is $\sigma$-continuous and 
$$\theta\big((x+\alpha 1)(y+\beta 1)\big)
\ = \ \Xi(xy + \beta x + \alpha y) 
\ = \ \check \Phi(x+\alpha 1)\theta(y) + \theta(x)\beta,$$
for any $x,y\in B^\omega$ and $\alpha, \beta\in \BC$. 
This implies that $\theta$ is a $(\check \Phi, \omega)$-derivation. 
The hypothesis gives $f_0\in E^*$ such that $\Xi(x) = \theta(x) = \check \Psi(x)(f_0)$ ($x\in B^\omega$) and the first conclusion is established. 

To show the second conclusion, we assume on the contrary that $I:=\overline{(B^\omega)^2}^\sigma \subsetneq B^\omega$.
Let $B_0 := B/I$ and $\sigma_0$ be the quotient topology on $B_0$. 
Denote by $q_0:B \to B_0$ the quotient map. 
Consider any non-zero Banach space $E$, and define a $\sigma$-continuous anti-representation $\Psi: B\to \CL(E)$ by $\Psi(a)x := \omega(a)x$ ($a\in B;x \in E$). 
Since $I\subsetneq B^\omega$ and $B = B^\omega \oplus \BC 1$, we know that the dimension of $B_0$ is strictly greater than one. 
Thus, there is a non-zero $\sigma_0$-continuous linear map $\theta:B_0\to E^*$ satisfying $\theta(q_0(1)) = 0$. 
Furthermore, for any $a,b\in B$, there exist unique elements $x,y\in B^\omega$ with $a = x + \omega(a) 1$ and $b = y + \omega(b) 1$.  
As $q_0(xy) = 0$, we see that 
$$\theta(q_0(ab)) = \omega(a)\theta(q_0(y)) + \theta(q_0(x))\omega(b) 
= \check\Psi(a)\theta(q_0(b)) + \theta(q_0(a))\omega(b),$$ 
which means that $\theta\circ q_0$ is a $(\check \Psi, \omega)$-derivation. 
However, $\theta\circ q_0$ is not inner (as the only inner $(\check \Psi, \omega)$-derivation is zero), and we have a contradiction.  

\smnoind
(c) $(1)\Leftrightarrow (2)$. 
This equivalence follows from Lemma \ref{lem:char-amen-bai} (notice that $\ti B$ is unital) and the following general facts. 
Any non-degenerate $\ti \sigma$-continuous anti-representation of $\ti B$ is of the form $\ti \Phi$ for a $\sigma$-continuous anti-representation $\Phi$ of $B$. 
Moreover, for any $\sigma$-continuous anti-representations $\Phi$ and  $\Psi$ of $B$, the assignment 
$\theta \mapsto \ti\theta$, where 
$\ti\theta((a, \lambda)) := \theta(a)$, 
is a bijection from the set of $\sigma$-continuous $1$-cocycles for $(\Phi, \Psi)$ to the set of $\ti \sigma$-continuous $1$-cocycles  for $(\ti \Phi, \ti \Psi)$ and every  $\sigma$-continuous $(\Phi, \Psi)$-derivation is inner if and only if every $\ti \sigma$-continuous $(\ti \Phi, \ti \Psi)$-derivation is inner. 

\noindent
$(2)\Leftrightarrow (3)$. 
This follows from parts (a) and (b). 
\end{prf}

\medskip

If $B$ is a $C^*$-algebra and $\omega$ is a $^*$-character, then part (c) above and Lemma \ref{lem:two-elem-inner}(b) (it is well-known that every $C^*$-algebra has a contractive approximate identity) tells us that $B$ is $\omega$-amenable. 

\medskip

Furthermore, parts (b) and (c) of the above, together with Lemma \ref{lem:two-elem-inner}(b), tell us that the $\omega$-amenability of $B_\sigma$ stands between $\ti B^{\ti \omega}$ having a $\ti \sigma$-bounded right $\ti \sigma$-approximate identity and $(\ti B^{\ti \omega})^2$ being $\ti \sigma$-dense in $\ti B^{\ti \omega}$. 
Moreover, if $B$ is unital, then parts (a) and (b) 
above imply that $B^\omega_{\sigma}$ has property (CB) if and only if $B_\sigma$ is $\omega$-amenable.

\medskip

On the other hand, Lemmas \ref{lem:bai-B-psi} and \ref{lem:two-elem-inner}(b) as well as Proposition \ref{prop:CB<->char-amen}(a) give the following partial generalization of \cite[Proposition 2.1]{KLP08}. 

\medskip

\begin{cor}\label{cor:CB-char-amen}
If $\ker \omega$ has a $\sigma$-bounded $\sigma$-approximate identity, then $B_\sigma$ is $\omega$-amenable. 
\end{cor}

\medskip

Let us recall that a Banach algebra $D$ is a \emph{$F$-algebra} (or a \emph{Lau algebra}) if there exists a von Neumann algebra structure on the dual space $D^*$ such that the identity $\varphi\in D^*$ (with respect to the von Neumann algebra structure) is a character on $D$ (see \cite{Lau83}). 
In this case, $D$ is said to be \emph{left amenable} if it is $\varphi$-amenable in the sense of Definition \ref{defn:char-amen} (see \cite[p.167]{Lau83}). 
If, in addition, $D$ is a Banach $^*$-algebra and $\varphi$ is a $^*$-homomorphism, then we call $D$ a \emph{$F^*$-algebra}. 

\medskip

In the following, we denote by $\ti L^1(G)$ the unitalization of the Banach $^*$-algebra $L^1(G)$. 
Notice that the ideals $\ti L^1(G)^{\ti \varepsilon_G}$ and $L^1(G)^{\varepsilon_G}$ coincide with $\ti L^1_0(G)$ and $L^1_0(G)$ as in \eqref{eqt:defn-L^1_0} and \eqref{eqt:defn-ti-L^1_0}, respectively. 

\medskip

\begin{thm}\label{thm:CB<->amen-group}
The following are equivalent for a locally compact group $G$. 

\begin{enumerate}
\item $G$ is amenable. 
\item The Banach algebra $\ti L^1_0(G)$ 
has property $(CB)$. 
\item The Banach algebra $L_0^1(G)$ has property $(CB)$. 
\item For any norm-continuous anti-representation $\mu:G\to \CL(E)$ 
and any norm continuous $L^1_0(G)$-module map $\Xi: L^1_0(G)\to E^*$, there exists $f\in E^*$ such that $\Xi(a)(x) = f(\hat{\mu}(a)x)$ ($a\in L_0^1(G); x\in E$). 
\end{enumerate}
\end{thm}
\begin{prf}
$(1) \Leftrightarrow (2)$. 
By \cite[Theorem 4.1]{Lau83}, $G$ is amenable if and only if $L^1(G)$ is left amenable, or equivalently, 
$\varepsilon_G$-amenable.
Now, the conclusion follows from Proposition \ref{prop:CB<->char-amen}(c). 

\smnoind
$(1) \Rightarrow (3)$. 
Since $G$ is amenable, \cite[Theorem 4.10]{Lau83} tells us that $L^1_0(G)$ has a norm-bounded right approximate identity, and Lemma \ref{lem:two-elem-inner}(b) implies that $L^1_0(G)$ has property $(CB)$. 

\smnoind
$(3) \Rightarrow (4)$. 
This is clear (by considering the restriction of the 
anti-representation $\hat \mu: L^1(G)\to \CL(E)$ to $L^1_0(G)$). 

\smnoind
$(4) \Rightarrow (1)$. 
Suppose that $\Phi:L^1(G) \to \CL(E)$ is a bounded non-degenerate anti-representation and $\theta$ is a $\|\cdot\|_{L^1(G)}$-continuous $(\check\Phi, \varepsilon_G)$-derivation.
By Lemma \ref{lem:cont-rep}(a), $\Phi$ is defined by a norm-continuous anti-representation $\mu:G\to \CL(E)$. 
Moreover, the hypothesis  and the argument for Proposition \ref{prop:CB<->char-amen}(a) tell us that $\theta$ is inner. 
Thus, Lemma \ref{lem:char-amen-bai} implies that $L^1(G)$ is $\varepsilon_G$-amenable, and hence $G$ is amenable (by \cite[Theorem 4.1]{Lau83}). 
\end{prf}

\medskip

The following result follows from similar arguments as that of Proposition \ref{prop:CB<->char-amen}, except that 
we employ Lemma \ref{lem:sep-n.d.}(b) instead of Lemma \ref{lem:char-amen-bai}. 

\medskip

\begin{lem}\label{lem:FH-unital}
Let $A^\varepsilon := \ker\varepsilon$. 

\smnoind
(a) If $A$ has a $\tau$-bounded $\tau$-approximate identity and $A^\varepsilon_{\tau}$ has property $(CH)$, then $(A_\tau, \varepsilon)$ has property $(FH)$. 

\smnoind
(b) If $A$ is unital and $(A_\tau, \varepsilon)$ has property $(FH)$,  then $A^\varepsilon_{\tau}$ has property $(CH)$ and $(A^\varepsilon)^2$ is $\tau$-dense in $A^\varepsilon$. 

\smnoind
(c) The following statements are equivalent. 
\begin{enumerate}
\item $(A_\tau, \omega)$ has property $(FH)$. 

\item $(\ti A_{\ti \tau}, \ti \omega)$ has property $(FH)$. 

\item $\ti A^{\ti \varepsilon}_{\ti \tau}$ has property $(CH)$. 
\end{enumerate}
\end{lem}

\medskip

\begin{eg}\label{eg:prop-FH}
Let $(A_\tau, \varepsilon)$ and $(B_\sigma,\omega)$ be as in the 
beginning of this section such that $A$ is commutative. 
Let $E$ be a non-zero Banach space. 

\smnoind
(a) Suppose that $E$ has an anti-linear isometry $^*$ 
satisfying $(x^*)^* = x$ ($x\in E$).
We may equip $E$  with the zero product, turning it into a commutative Banach algebra. 
If $\varepsilon_0:\ti E\to \BC$ is the unital $^*$-homomorphism that vanishes on $E$, then Lemma \ref{lem:FH-unital}(b) tells us that $(\ti E, \varepsilon_0)$ does not have property $(FH)$.

\smnoind
(b) Suppose that $\nu:A\to \CL(\KH)$ is a $^*$-representation with $\ker \nu \nsubseteq \ker \varepsilon$, and $\theta$ is a $(\nu, \varepsilon)$-derivation.  
Fix any element $x\in \ker \nu$ with $\varepsilon(x) = 1$.
Since $A$ is commutative, one has $\theta(ax) = \theta(xa)$ and $\theta(a) = \theta(x) \varepsilon(a) - \nu(a) \theta(x)$.
Thus, $\theta$ is inner. 

\smnoind
(c) Suppose that $\Phi_\omega:B\to \CL(E)$ is the anti-representation given by $\Phi_\omega(a)(x) := \omega(a)x$ ($a\in B; x\in E$).
If $\overline{(B^\omega)^2}^{\sigma} = B^\omega$, then every $\sigma$-continuous $(\check \Phi_\omega, \omega)$-derivation $\theta$ is zero. 

In fact, choose any $u\in B$ with $\omega(u) = 1$. 
As $\theta ((B^\omega)^2) = \{0\}$ and $u^2-u\in B^\omega$, we know from $\overline{(B^\omega)^2}^{\sigma} = B^\omega$ that $\theta(u^2) = \theta(u)$.
On the other hand, as $\theta(u^2) = 2\theta(u)$, we conclude that $\theta(u) = 0$. 
Thus, $\theta = 0$ (because $B=B^\omega + \BC u$). 

\smnoind
(d) Let $\Phi_\omega$ be as in part (c). 
If there exists $p\in B$ 
satisfying $p^2=p$, $\omega(p) =1$ and $pa=\omega(a)p = ap$ ($a\in B$), then every $(\Phi_\omega, \omega)$-derivation $\theta$ is zero. 

In fact, the equality $\theta(p) = \theta(p^2) = 2\theta(p)$ tells us that $\theta(p) = 0$. 
Moreover, for any $a\in B$, if we set $a_0:= a - \omega(a)p \in B^\omega$, then
the relation $\theta((B^\omega)^2) = \{0\}$ implies
$$2\omega(a)\theta(a) 
 =  \theta(a^2)
 =  \theta\big((a_0+\omega(a)p)(a_0+ \omega(a)p)\big)
 = 0.$$
Thus, for any $c\in B^\omega$, we have $\theta(c) = \theta(c+p) = 0$ (as $\omega(c+p) = 1$) and we conclude that $\theta \equiv 0$. 

In particular, if $G$ is a compact group, then the only $(\Phi_{\varepsilon_G}, \varepsilon_G)$-derivation on $C_c(G)$ is zero (see the argument of Theorem \ref{thm:char-amen-cpt} in Section 5 below). 

\smnoind
(e) If the commutative algebra $A$ satisfies both $\overline{(A^\varepsilon)^2}^{\tau} = A^\varepsilon$ and the following ``$T$-like condition'': 
\begin{quotation}
any $\tau$-continuous $^*$-representation $\Phi:A\to \CL(\KK)$ satisfying $\ker \Phi \subseteq \ker \varepsilon$ actually contains $\varepsilon$, 
\end{quotation}
then $(A_\tau,\varepsilon)$ has property $(FH)$. 

In fact, let $\Phi:A\to \CL(\KH)$ be a $\tau$-continuous $^*$-representation and 
$$\KH_{\Phi, \varepsilon}:=\{\xi\in \KH: \Phi(a)\xi = \varepsilon(a)\xi, \text{ for any }a\in A\}.$$ 
If $\KK := \KH_{\Phi, \varepsilon}^\bot$ and $\nu:A\to \CL(\KK)$ is the $^*$-representation induced by $\Phi$, then one has $\ker\nu \nsubseteq \ker \varepsilon$ because of the property displayed above. 
The conclusion now follows from parts (b) and (c). 

Observe that this $T$-like condition is not a necessity for property $(FH)$, e.g., if $G$ is a locally compact group and $\delta_e\in C_0(G)^*$ is the evaluation at the identity $e\in G$, then 
$(C_0(G), \delta_e)$ has property $(FH)$ (by Corollary \ref{cor:CB-char-amen}). 
However, when $G$ is abelian, the above $T$-like condition will imply that $G$ is discrete. 

\smnoind
(f) If $A$ is a commutative $F^*$-algebra, $\tau$ is the norm topology and $\varepsilon$ is the identity of a von Neumann algebra structure on $A^*$, then 
Example (1) in \cite[p.168]{Lau83} tells us that $A$ is $\varepsilon$-amenable.
In particular, $(A_\tau, \varepsilon)$ has property $(FH)$. 

\smnoind
(g) By part (f), we know that $(A(G), \delta_e)$ has property $(FH)$, where $\delta_e$ is the evaluation at the identity $e\in G$. 
This shows that the assumption in Corollary \ref{cor:CB-char-amen} is not an absolute necessity since the existence of a bounded approximate identity in $A(G)^{\delta_e}$ is equivalent to $G$ being amenable (see \cite[Theorem 4.10]{Lau83}). 
\end{eg}

\medskip

One may generalize Example \ref{eg:prop-FH}(g) to the case of amenable locally compact quantum groups. 
The definition, notations and properties of a locally compact quantum group $\BG$ can be found in, e.g., \cite{HNR}, \cite{KV1} and \cite{Tim}. 
Since the dual space of the Banach algebra $L^1(\BG)$ is the von Neumann algebra $L^\infty(\BG)$ and the identity of $L^\infty(\BG)$ is a homomorphism on $L^1(\BG)$, we know that $L^1(\BG)$ is a $F$-algebra. 
Recall also that $\BG$ is said to be \emph{amenable} if $L^\infty(\BG)$ has an ``invariant mean'' (the readers may consult, e.g., \cite{BT} for the precise meaning of invariant mean and the properties of amenable quantum groups). 
Thus, by \cite[Theorem 4.1]{Lau83}, $L^1(\BG)$ is left amenable when $\BG$ is amenable. 
Moreover, $\BG$ is said to be \emph{of Kac type} if its antipode is bounded and its modular element is affiliated with the center of $L^\infty(\BG)$.  
In this case, the bounded antipode turns $L^1(\BG)$ into a $F^*$-algebra. 
These give part (a) of the following result.
Note that part (b) follows from \cite[Theorem 5.1]{Kyed} as well as parts (a) and (b) of Lemma \ref{lem:FH-unital}. 

\medskip

\begin{cor}\label{cor:lcqg}
Let $\BG$ be a locally compact quantum group and $\varepsilon_\BG$ be the trivial one-dimensional representation of $C_0^{\rm u}(\widehat\BG)$.

\smnoind
(a) If $\BG$ is amenable and of Kac type, then $(L^1(\BG), \varepsilon_\BG)$ has property $(FH)$.

\smnoind
(b) Suppose that $\BG$ is discrete and separable. 
Then $\BG$ has property $(T)$ if and only if the $^*$-algebra 
\begin{equation*}\label{eqt:defn-Pol^0}
{\rm Pol}^0(\widehat{\BG}):=\big\{x\in {\rm Pol}(\widehat{\BG}): \varepsilon_\BG(x)=0\big\}, 
\end{equation*}
when equipped with the discrete topology, has property $(CH)$. 
\end{cor}

\medskip

Observe, however, that the amenability assumption in part (a) above is not an absolute necessity (see e.g.\ Theorem \ref{thm:FH=T}(a) below). 

\medskip

Now, we go back to the consideration of the locally compact group $G$.
Obviously, $C_c(G)$ is the vector space inductive limit of the system $\{C_K(G)\}_{K\in \KK(G)}$, and we consider $\IT$ to be the locally convex inductive topology on $C_c(G)$, when all $C_K(G)$ are equipped with the $L^1$-norms (see e.g. \cite[\S II.6]{Sch}).
Observe that $\IT$ is strictly finer than the $L^1$-norm on $C_c(G)$ and a linear map from $C_c(G)$ to a Banach space is $\IT$-continuous if and only if it is locally bounded.
Moreover, for a fixed $f\in C_c(G)$, the maps $g\mapsto f\ast g$, $g\mapsto g\ast f$ and $g\mapsto g^*$ are $\IT$-$\IT$-continuous.
Hence, $C_c(G)_\IT$ is a locally convex $^*$-algebra. 

\medskip

\begin{lem}\label{lem:T-bdd-T-app-id}
Let $B\subseteq C_c(G)$ be a subalgebra and $K\subseteq G$ be a compact subset. 
If $\{h_i\}_{i\in \KI}$ is a $\|\cdot\|_{L^1(G)}$-bounded $\|\cdot\|_{L^1(G)}$-approximate identity in $B$ such that $\supp h_i\subseteq K$ ($i\in\KI$), then $\{h_i\}_{i\in \KI}$ is a $\IT$-bounded $\IT$-approximate identity in $B$. 
\end{lem}
\begin{prf}
Consider any $f\in B$ and 
set $L:= K\cdot (\supp f)\cup \supp f \cup K$.
For any $\IT$-neighbourhood $V\subseteq C_c(G)$ of zero, there exists $\delta > 0$ such that 
$$\{g\in C_L(G): \|g\|_{L^1(G)} \leq \delta\} \subseteq V.$$ 
Thus, $\{h_i\}_{i\in \KI}$ is $\IT$-bounded. 
Moreover, as $h_i * f - f\in C_L(G)$ for any $i\in \KI$, we know that $h_i * f - f\in V$ when $i$ is large enough, and $\{h_i\}_{i\in \KI}$ is a left $\IT$-approximate identity. 
In a similar fashion, one can show that $\{h_i\}_{i\in \KI}$ is a right $\IT$-approximate identity.
\end{prf}

\medskip

Consequently, if $U$ is a fixed element in $\CU$, then the net $\{h_V\}_{V\in \CU(U)}$ as in Section 2 is a $\IT$-bounded $\IT$-approximate identity of $C_c(G)$.

\medskip
 
As in the above, we denote by $\ti C_c(G)$ the unitalization of $C_c(G)$. 
Note that the ideals $\ti C_c(G)^{\ti \varepsilon_G}$ and $C_c(G)^{\varepsilon_G}$ (see Lemma \ref{lem:FH-unital}) coincide with $\ti C_c^0(G)$ and $C_c^0(G)$ as in 
\eqref{eqt:defn-ti-C_c^0} and \eqref{eqt:defn-C_c^0}, respectively. 

\medskip

\begin{thm}\label{thm:FH=T}
Let $G$ be a locally compact group.

\smnoind
(a) $(L^1(G), \varepsilon_G)$ has property $(FH)$.

\smnoind
(b) In the case when $G$ is second countable, the following statements are equivalent. 
\begin{enumerate}
\item $G$ has property $(T)$. 
\item $(C_c(G)_\IT, \varepsilon_G)$ has property $(FH)$.
\item $\ti C_c^0(G)_{\ti \IT}$ has property $(CH)$.
\item For any norm-continuous unitary representation $\pi:G\to \CL(\KH)$, if $\Xi:\ti C_c^0(G)\to \KH$ is a $\ti \IT$-continuous $\ti C_c^0(G)$-module map, then there is $\xi_0\in \KH$ such that 
$\Xi(a) = \hat \pi(a)\xi_0$ ($a\in \ti C_c^0(G)$). 
\end{enumerate}

\smnoind
(c) Suppose that $G$ is second countable. 
If $C_c^0(G)_\IT$ has property $(CH)$, then $G$ has property $(T)$. 
\end{thm}
\begin{prf}
(a) This follows from Lemma \ref{lem:cf-1-cocyc}(d) and \cite[Proposition 2.2.9]{BHV}. 

\smnoind
(b) $(1) \Leftrightarrow (2)$. 
This follows from Lemmas \ref{lem:sep-n.d.}(b) and \ref{lem:cont-rep} as well as Theorem  \ref{thm:alter-prop-T-1}. 

\smnoind
$(2)\Leftrightarrow (3)$. 
This follows from Lemma \ref{lem:FH-unital}(c). 

\smnoind
$(3)\Rightarrow (4)$. 
This is clear. 

\smnoind
$(4)\Rightarrow (1)$. 
By the equivalence of Statements (1) and (2), it suffices to show that $(C_c(G)_\IT, \varepsilon_G)$ has property $(FH)$. 
Let $A:= C_c(G)$. 
Suppose that $\Phi: A\to \CL(\KH)$ is a non-degenerate $\IT$-continuous $^*$-representation and $\theta:A\to \KH$ is a $\IT$-continuous $(\Phi, \varepsilon_G)$-derivation. 
By Lemma \ref{lem:cont-rep}(b), there is a norm-continuous unitary representation $\pi:G\to \CL(H)$ with $\Phi = \hat \pi$. 
As in the proof of Proposition \ref{prop:CB<->char-amen}(c), if $\ti \theta:\ti A\to \KH$ is defined by 
$$\ti \theta(a+\alpha 1) := \theta(a) \qquad (a\in A; \alpha\in \BC),$$ 
then $\ti \theta$ is a $(\ti \Phi, \ti \varepsilon_G)$-derivation. 
Set $\Psi: = \ti \Phi|_{\ti A^{\ti \varepsilon_G}}$ and $\Xi:=\ti \theta|_{\ti A^{\ti \varepsilon_G}}$. 
Then $\Xi$ is a $\ti \IT$-continuous $\ti A^{\ti \varepsilon_G}$-module map, and one obtains $\xi_0\in\KH$ satisfying the relation in Statement (4). 
By the argument of Proposition \ref{prop:CB<->char-amen}(a), we see that $\ti \theta$ is inner and so $\theta$ is inner. 
Now, the conclusion follows from Lemma \ref{lem:sep-n.d.}(b). 

\smnoind
(c) This is consequence of part (b) and Lemma \ref{lem:FH-unital}(a). 
\end{prf}

\medskip

Part (a) above tells us that the cohomology characterization for property $(T)$ as in part (b) cannot be done at the group Banach algebra level (and hence neither at the full group $C^*$-algebra level; see also the paragraph following Proposition \ref{prop:CB<->char-amen}).

\medskip

We do not know whether $G$ having property $(T)$ will imply $C_c^0(G)_\IT$ has property $(CH)$. 
Note that the corresponding equivalence in Theorem \ref{thm:CB<->amen-group} depends on \cite[Theorem 4.10]{Lau83}, 
which ensures the existence of a bounded approximate identity. 
However, the corresponding fact is not true for property $(FH)$. 
In fact, if $G$ is a second countable locally compact group, then the existence of a 
$\IT$-bounded $\IT$-approximate identity in $C_c^0(G)$ will imply that $C_c(G)_\IT$ is $\varepsilon_G$-amenable (see Corollary \ref{cor:CB-char-amen}) and hence $G$ is compact, because of Theorem \ref{thm:char-amen-cpt}(b) in the next section. 
	
\bigskip

\section{An application to fixed point property for affine actions}

\medskip

In this section, we will use the ideas and arguments 
in the previous sections to obtain some fixed point results. 
Let us first set some notations. 

\medskip

Suppose that $X$ is a (complex) vector space.
Recall that a map $S:X\to X$ is \emph{affine-linear} if there exist a (complex) linear map $S^\Kl:X\to X$ as well as an element $x_S\in X$ such that $S(x) = S^\Kl(x) + x_S$ ($x\in X$). 
Note that when $S$ is bijective, $S^{-1}$ is also affine-linear. 
On the other hand, for a convex subset $C$ of a vector space, a map $\Lambda: C\to C$ is
said to be \emph{affine} if $\Lambda\big( t x + (1-t)y\big)  = t \Lambda(x) + (1-t) \Lambda(y)$ ($x,y\in C; t\in [0,1]$).

\medskip

We recall Day's fixed point theorem as follows (see \cite{Day}, \cite{Day-Cor} and \cite[p.49]{Green}): 
\begin{quote}
$G$ is amenable if and only if any continuous affine action of $G$ on a non-empty compact convex subset $K$ of a locally convex space has a fixed point.
\end{quote}

\medskip

The following proposition can be regarded as is a variant of Day's fixed point theorem of amenable groups concerning affine-linear actions  on a dual Banach space $E^*$ rather than affine actions on weak$^*$-compact convex subsets sets of $E^*$. 
It should be noted that ``affine actions'' in 
this result cannot be replaced by ``linear actions'' since any linear action always has a common fixed point, namely ``0''. 

\medskip

\begin{prop}\label{prop:amen}
$G$ is amenable if and only if any weak$^*$-continuous affine-linear action $\alpha$ of $G$ on any dual Banach space $E^*$ with one norm-bounded orbit (and equivalently, with all orbits being bounded) has a fixed point. 
\end{prop}
\begin{prf}
$\Rightarrow)$. 
Consider 
an element $\varphi_0\in E^*$ with the orbit $O:=\{\alpha_t(\varphi_0):t\in G\}$ being norm-bounded. 
Set $C$ to be the weak$^*$-closure of the convex hull of $O$. 
As $O$ is norm-bounded, $C$ is weak$^*$-compact. 
Now, Day's fixed point theorem produces a fixed point as required.  

\smnoind
$\Leftarrow)$. 
Let $\Phi:L^1(G)\to \CL(E)$ be a non-degenerate bounded anti-representation, $\check \Phi:L^1(G)\to \CL(E^*)$ be the induced map and $\theta:L^1(G)\to E^*$ be a bounded $(\check \Phi, \varepsilon_G)$-derivation. 
By Lemma \ref{lem:cont-rep}(a), 
$\Phi = \hat \mu$ for a norm-continuous anti-representation $\mu:G\to \CL(E)$ and we set $\pi_t := \mu_t^*\in \CL(E^*)$ ($t\in G$).

On the other hand, Lemma \ref{lem:cf-1-cocyc}(b)  allows us to define an affine-linear action $\gamma_t(\varphi):= \pi_t(\varphi) + \D(\theta)(t)$ ($\varphi\in E^*;t\in G$). 
As $\pi$ is a weak$^*$-continuous action on $E^*$ and $\D(\theta)$ is weak-$^*$-continuous, we know that the affine-linear action $\gamma$ is weak$^*$-continuous.

Furthermore, the boundedness of $\theta$ and Relation \eqref{eqt:defn-D-theta} give 
$$\|\D(\theta)(r)(x)\|\ \leq\ \|\theta\| \|x\| \qquad (x\in E; r\in G),$$ 
and $\{\D(\theta)(t):t\in G\}$ is a bounded subset of $E^*$.  
Since the subset $\{\|\pi_t\|: t\in G\}$ is also bounded (see e.g.\ Relation \eqref{eqt:loc-bdd-pi}), we conclude that all orbits of $\gamma$ are norm-bounded. 

The hypothesis now produces a fixed point $\varphi_0\in E^*$ for $\gamma$, and it is not hard to check that $\theta(f) = \varepsilon_G(f)\varphi_0 - \check \Phi(f)\varphi_0$. 
By Lemma \ref{lem:char-amen-bai}, the Banach algebra $L^1(G)$ is left amenable and hence $G$ is amenable (by \cite[Theorem 4.1]{Lau83}). 
\end{prf}

\medskip

It was recently shown in \cite{GM} that a locally compact $\sigma$-compact group $G$ is amenable if and only if every continuous affine-linear action of $G$ on a separable real Hilbert space with a bounded orbit has a fixed point. 

\medskip

Recall that a locally convex space $X$ is \emph{quasi-complete} if every closed and bounded subset of $X$ is complete. 
Observe that in this case, the closure of any totally bounded subset of $X$ is compact. 
Examples of quasi-complete locally convex spaces include all Banach spaces with the norm-topologies and all dual Banach spaces with the weak$^*$-topologies. 

\medskip

We also have the following corollary, which give an analogue of Day's fixed point theorem for compact 
groups. 

\medskip

\begin{cor}\label{cor:fix-pt}
The following statements are equivalent for a $\sigma$-compact locally compact group $G$. 
\begin{enumerate}
\item $G$ is compact. 
\item Each continuous affine action $\gamma$ of $G$ on a non-empty closed convex subset $C$ of a quasi-complete locally convex space $X$ has a fixed point (in $C$). 
\item For every non-degenerate bounded anti-representation $\Phi:L^1(G)\to \CL(E)$, any $\IT$-continuous $(\check \Phi, \varepsilon_G)$-derivation $\theta: C_c(G)\to E^*$ is inner. 
\end{enumerate}
\end{cor}
\begin{prf}
$(1)\Rightarrow (2)$. 
Pick any $x\in C$.
As the set $\{\gamma_t(x): t\in G\}$ is totally bounded, its convex hull $C_0$ is also totally bounded  (see e.g.\ \cite[Theorem 3.24]{Rudin}), and hence $\overline{C_0}$ is compact. 
Now, Day's fixed point theorem produces a $\gamma$-fixed point in $\overline{C_0}$. 

\smnoind
$(2)\Rightarrow (3)$. 
As in the argument of Proposition \ref{prop:amen}, the map $\theta$ induces a weak$^*$-continuous action $\gamma$ of $G$ on $E^*$ by affine-linear maps. 
Now, this implication follows from the argument of Proposition \ref{prop:amen} (by taking 
$C=E^*$). 

\smnoind
$(3)\Rightarrow (1)$. 
Let $\pi$ be a norm-continuous unitary rerpesentation of $G$ on a Hilbert space $\KH$. 
By considering $\Phi(a) := \hat \pi(a)^*\in \CL(\KH^*)$ ($a\in L^1(G)$), we conclude from Statement (3) that $H^1(G;\pi) = (0)$.
Thus, the Delorme-Guichardet theorem implies that $G$ has property $(T)$. 
On the other hand, the argument of Proposition \ref{prop:amen} also implies that $G$ is amenable.
Consequently, $G$ is compact. 
\end{prf}

\medskip


\begin{thm}\label{thm:char-amen-cpt}
Let $G$ be a locally compact group and $C_c^0(G)$ is as in \eqref{eqt:defn-C_c^0}.
Consider the following statements. 
\begin{enumerate}
\item $G$ is compact.
\item There exist a compact subset $K\subseteq G$ and a $\|\cdot\|_{L^1(G)}$-bounded $\|\cdot\|_{L^1(G)}$-approximate identity $\{h_i\}_{i\in \KI}$ in  $C_c^0(G)$ with $\supp h_i \subseteq K$ ($i\in \KI$). 
\item $C_c^0(G)$ has a $\IT$-bounded $\IT$-approximate identity. 
\item $C_c^0(G)_\IT$ has property $(CB)$. 
\item $C_c(G)_\IT$ is $\varepsilon_G$-amenable. 
\end{enumerate}
(a) One has $(1)\Rightarrow (2) \Rightarrow (3) \Rightarrow (4) \Rightarrow (5)$. 

\smnoind
(b) If $G$ is $\sigma$-compact, then $(5)\Rightarrow (1)$.
\end{thm}
\begin{prf}
(a) $(1) \Rightarrow (2)$. 
Let $1_G\in C_c(G)$ be the constant one function, and $\{h_V\}_{V\in \CU}$ be the 
$\|\cdot\|_{L^1(G)}$-bounded $\|\cdot\|_{L^1(G)}$-approximate identity in $C_c(G)$ as in Section 2.
Then $1_G * f = \varepsilon_G(f)1_G = f* 1_G$ ($f\in C_c(G)$). 
Thus, if we set  $g_V := h_V - 1_G\in C_c^0(G)$ ($V\in \CU)$, 
then $\|g * g_V - g\|_{L^1(G)} \to 0$ and $\|g_V * g - g\|_{L^1(G)} \to 0$ ($g\in C_c^0(G)$). 

\noindent
$(2) \Rightarrow (3)$. 
This implication follows from Lemma \ref{lem:T-bdd-T-app-id}. 

\noindent
$(3) \Rightarrow (4)$. 
This follows from Lemma \ref{lem:two-elem-inner}(b). 

\noindent
$(4) \Rightarrow (5)$. 
This follows from Proposition \ref{prop:CB<->char-amen}(a) (since 
$C_c(G)$ has a $\IT$-bounded $\IT$-approximate identity by Lemma \ref{lem:T-bdd-T-app-id}). 

\smnoind
(b) This part follows from Lemma \ref{lem:char-amen-bai} and the argument for Corollary \ref{cor:fix-pt}. 
\end{prf}

\bigskip

\section{Open questions}

\medskip

The following is one of the motivating questions of this work. 

\medskip

\begin{quest}
Can one extend the Delorme-Guichardet theorem to locally compact quantum groups? 
\end{quest}

\medskip

Before answering this question, one needs to consider the following.

\medskip

\begin{quest}
Is it possible to define an analogue of $C_c(G)$ for a locally compact quantum group $\BG$? 
\end{quest}

\medskip

In fact, we do not know the answer for this question even in the case when $\BG$ is the dual quantum group of a locally compact group. 
More precisely, in the case of a locally compact group $G$, 
we do not know how to define a canonical dense $^*$-subalgebra $A_0(G)$ of $A(G)$ such that $A_0(G) \cong C_c(\widehat{G})$ under the Fourier transform, when $G$ is abelian. 

\medskip

In the following, we list some more questions that are related to the results in the paper. 
A Banach right $G$-module $E$ is said to be \emph{contractive} if $\|x\cdot t\| \leq \|x\|$ ($x\in E;t\in G$). 
Note that when $\pi:G\to \CL(\KH)$ is a norm-continuous unitary representation, the induced Banach right $G$-module structure on $\KH^*$ is always contractive. 
In this respect, the following is a natural question arising from Theorem \ref{thm:FH=T}(b) and Theorem \ref{thm:char-amen-cpt}. 

\medskip

\begin{quest}
Suppose that $G$ is second countable. 
Can one describe the topological or analytical property for $G$ under which for any continuous contractive Banach right $G$-module $E$, any weak$^*$-continuous $1$-cocycle $c: G \to E^*$ is a $1$-coboundary?
\end{quest}

\medskip

The above property is closely related to Property $(F_B)$, as introduced in \cite[Definition 1.2]{BFGM}, for a Banach space $B$. 
However, apart from the obvious difference that one considers all Banach spaces instead of a fixed Banach space (as in the case of Property $(F_B)$), the above property concerns with cocycle taking value in the dual Banach spaces instead of the original Banach space. 

\medskip

Another natural question is whether one can remove the $\sigma$-compactness from Corollary \ref{cor:fix-pt} (and hence Theorem \ref{thm:char-amen-cpt} as well). 

\medskip

\begin{quest}\label{quest:rem-sig-cpt}
Can one remove the $\sigma$-compact assumption in Corollary \ref{cor:fix-pt}?
\end{quest}

\medskip

A related question is the following. 

\medskip

\begin{quest}
If $G$ is amenable and has property (FH), will $G$ be compact?
\end{quest}

\medskip

The following question is also interesting. 

\medskip

\begin{quest}
Does the converse of Theorem \ref{thm:FH=T}(c) hold? 
\end{quest}

\medskip

Finally, one may also consider the following question. 

\begin{quest}
Can the left amenability of a semigroup be characterized by affine actions on dual Banach spaces as in Proposition \ref{prop:amen} (see also \cite{LZ-08} and  
\cite{LZ-12})?
\end{quest}

\bigskip

\appendix \section {Some known facts}

\medskip

This appendix contains three probably well-known results.
Since we do not find them explicitly stated in the literature, we give their complete arguments here for the benefit of the reader. 
The first two results concern with a locally compact group $G$, and the third one 
concerns with the analogues of two well-known facts in Banach algebras. 

\medskip

\begin{lem}\label{lem:equi-cont}
	If $K$ is a compact subset of $G$ and $f\in L^1(G)$,
	then $$\sup_{r\in K} \|\rho_r(f)\ast h_V - \rho_r(f)\|_{L^1(G)} \to 0 \quad (\text{along } V\in \CU).$$
\end{lem}
\begin{prf}
	Note that $C:=\{\rho_r(f): r\in K\}$ is a $\|\cdot\|_{L^1(G)}$-compact subset of $L^1(G)$. 
	Given $\epsilon >0$, there exists a finite subset $F\subseteq C$ with ${\rm dist}(g, F) < \epsilon$ whenever $g\in C$. 
	For this finite set $F$, one can find $V_0\in \CU$ such that for any $V\subseteq V_0$, one has $\|h\ast h_V - h\|_{L^1(G)} < \epsilon$ ($h\in F$),  and hence $\|g\ast h_V - g\|_{L^1(G)} < 3\epsilon$ ($g\in C$). 
\end{prf}

\medskip

\begin{lem}\label{lem:cont-rep}
	(a) If $E$ is a Banach space, and $\Psi:C_c(G)\to \CL(E)$ is a non-degenerate locally bounded representation (respectively, anti-representation), then there exists a unique norm-continuous representation (respectively, anti-representation) $\mu: G\to \CL(E)$ such that $\Psi = \hat\mu$.
	
	\smnoind
	(b) If $\KH$ is a Hilbert space, and $\Phi:C_c(G)\to \CL(\KH)$ is a non-degenerate locally bounded $^*$-representation, then there exists a unique norm-continuous unitary representation $\pi$ of $G$ on $\KH$ with $\Phi = \hat \pi$.
\end{lem}
\begin{prf}
	(a) Fix a $U\in \CU$ and consider the approximate identity $\{h_V\}_{V\in \CU(U)}$ as in Section 2.
	Let $s\in G$ and let $W$ be an open neighborhood of $s$ with compact closure.
	Suppose that $f_1,...,f_n\in C_c(G)$ and $x_1,...,x_n\in E$.
	If $L$ is the compact set ${\bigcup}_{i=1}^n \overline{W}\cdot \overline{U}\cdot (\supp f_n)$, then the local boundedness of $\Psi$ produces $\kappa>0$ with $\|\Psi(g)\| \leq \kappa \|g\|_{L^1(G)}$ ($g\in C_L(G)$).
	For any $r\in W$ and $V\in \CU(U)$, we have
	$$\Big\|\Psi\big(\lambda_r(h_V)\big)\Big({\sum}_{i=1}^n \Psi(f_i)x_i\Big) - {\sum}_{i=1}^n \Psi\big(\lambda_r(f_i)\big)x_i\Big\|
	\ \leq \ \kappa {\sum}_{i=1}^n \| h_V*f_i - f_i\|_{L^1(G)}\|x_i\|,$$
	which converges to zero uniformly for all $r\in W$.
	Moreover, as $\|\Psi(\lambda_r(h_V))\|\leq \kappa$ whenever $r\in W$ and $V\in \CU(U)$, we know that for any $x\in E$, the family of nets $\big\{\Psi\big(\lambda_r(h_V)\big)x\big\}_{V\in \CU(U)}$  is uniformly norm-Cauchy for all $r\in W$ (as $\Psi$ is non-degenerate).
	Thus, $\big\{\Psi\big(\lambda_r(h_V)\big)x\big\}_{V\in \CU(U)}$ norm-converges to an element $\mu(r)x\in E$ uniformly for $r\in W$. 
	Since 
	\begin{equation}\label{eqt:loc-bdd-pi}
	\|\mu(r)x\|\ \leq\ \kappa\|x\| \qquad (r\in W; x\in E),  
	\end{equation}
	we know that $\mu(r)\in \CL(E)$. 
	Furthermore, as 
	\begin{equation}\label{eqt:rep}
	\mu(r)\Big({\sum}_{i=1}^n \Psi(f_i)x_i\Big) = {\sum}_{i=1}^n \Psi\big(\lambda_r(f_i)\big)x_i,
	\end{equation}
	the norm-density of $\Psi(C_c(G))E$ in $E$ implies that $\mu(r)x$ does not depend on the choices of $U$, $\{h_V\}_{V\in \CU(U)}$ nor $W$ (so long as $W$ contains $r$).
	Since $r\mapsto \Psi\big(\lambda_r(h_V)\big)x$ are continuous maps from $W$ to $E$, for any $V\in \CU(U)$, and the convergence to $\mu(r)x$ is uniform for all $r\in W$, we see that $\mu$ is norm-continuous (as $s$ and $W$ are arbitrary).
	Moreover, Equality \eqref{eqt:rep} and the non-degeneracy of $\Psi$ also tell us that $\mu$ is a representation of $G$ with $\hat \mu = \Psi$. 
	
	Suppose that $\nu:G\to \CL(E)$ is another norm-continuous representation of $G$ with $\hat \nu = \Psi$. 
	For $f,g\in C_c(G)$ and $x\in E$, one has 
	$$\int_G g(r)\nu(r)\Psi(f)x\ \! dr\ =\ \Psi(g\ast f)x\ =\ \int_G g(r) \Psi(\lambda_r(f))x\ \!dr.$$
	This implies that $\nu(r)\Psi(f)x = \Psi(\lambda_r(f))x = \mu(r)\Psi(f)x$ for all $r\in G$ (because of the norm continuity of both $\nu$ and $\mu$). 
	Thus, $\nu = \mu$ 
	(again, thanks to the non-degeneracy of $\Psi$). 
	
	\smnoind
	(b) By part (a), one can find a norm-continuous representation $\pi:G\to \CL(\KH)$ such that $\Phi = \hat \pi$. 
	Equality \eqref{eqt:rep} tells us that for any $r\in G$, $f,g\in C_c(G)$ and $\eta, \zeta\in \KH$, one has, 
	\begin{eqnarray*}
		\big\la \pi(r^{-1})\Phi(f)\eta, \Phi(g)\zeta\big\ra 
		& = & \big\la \Phi(g^*)\Phi\big(\lambda_{r^{-1}}(f)\big)\eta, \zeta\big\ra 
		\ = \ \big\la \Phi\big((g^*\ast \lambda_{r^{-1}}(f)\big)\eta, \zeta\big\ra \\
		& = & \big\la \Phi\big((\lambda_{r}(g))^*\ast f\big)\eta, \zeta\big\ra
		\ = \ \big\la \Phi(f)\eta, \pi(r)\Phi(g)\zeta\big\ra, 
	\end{eqnarray*}
	which shows that $\pi(r^{-1}) = \pi(r)^*$ and hence $\pi$ is a unitary representation. 
\end{prf}

\medskip

\begin{lem}\label{lem:bai-B-psi}
Suppose that $B_\sigma$ is a locally convex algebra. 

\smnoind
(a) Let $\omega:B\to \BC$ be a $\sigma$-continuous character.  
If $\ker\omega$ has a $\sigma$-bounded left (respectively, right) $\sigma$-approximate identity $\{e_i\}_{i\in \KI}$, then $B$ has a $\sigma$-bounded left (respectively, right) $\sigma$-approximate identity. 

\smnoind
(b) If $B$ has both a $\sigma$-bounded left $\sigma$-approximate identity $\{a_i\}_{i\in \KI}$ and a $\sigma$-bounded right $\sigma$-approximate identity $\{b_i\}_{j\in \KJ}$, then $B$ has a $\sigma$-bounded $\sigma$-approximate identity. 
\end{lem}
\begin{prf}
(a) Pick any $u\in B$ with $\omega(u) = 1$. 
Since $x - ux\in \ker\omega$ 
(respectively, $x - xu\in \ker\omega$) for any $x\in B$, it is easy to check that 
$\{e_i - e_iu +u\}_{i\in \KI}$ (respectively, $\{e_i - ue_i +u\}_{i\in \KI}$) is a $\sigma$-bounded left (respectively, right) $\sigma$-approximate identity for $B$. 

\smnoind
(b) For any $x\in B$, if $V$ is a neighbourhood  of zero in $B$, then
there are neighbourhoods $W_1$, $W_2$ and $W_3$ of zero in $B$ satisfying $W_1 + W_2\cdot W_3\subseteq V$. 
The boundedness of $\{b_j\}_{j\in \KI}$ implies the existence of $N\in \BN$ with $\{b_j\}_{j\in \KJ}\subseteq N W_3$.
Moreover, 
there is $i_0\in \KI$ such that  whenever $i\geq i_0$, one has $x - xa_i \in W_1\cap \frac{1}{N}W_2$, and hence $x - x(a_i + b_j - a_ib_j)\in V$. 
Similarly, $\{a_i + b_j - a_ib_j\}_{(i,j)\in \KI\times \KJ}$ is also a left $\sigma$-approximate identity. 
\end{prf}

\bigskip

\section*{Acknowledgement}

\medskip

The first and the last named authors are supported by the National Natural Science Foundation of China (11071126 and 11471168).
The second named author is supported by NSERC grant MS100.

Parts of this work was done during the visit of the second and the last named authors to the Fields Institute during the Thematic Programme on Abstract Harmonic Analysis and Operator Algebras, 2014.

\end{document}